\documentclass[aop,seceqn,MSNbibl,citesort,dvips]{arximspdf}

\doi{10.1214/10-AOP577}
\volume{39}
\issue{3}
\pubyear{2011}
\firstpage{1027}
\lastpage{1060}

\makeatletter

\newcommand{\zz}{{\mathbb Z}}
\newcommand{\re}{\operatorname{Re}}
\newcommand{\im}{\operatorname{Im}}
\newcommand{\p}{{\mathbb P}}
\newcommand{\q}{{\mathbb Q}}
\newcommand{\e}{{\mathbb E}}
\newcommand{\rr}{{\mathbb R}}
\newcommand{\ccc}{{\mathbb C}}
\newcommand{\dd}{{ d}}
\newcommand{\ii}{{ i}}
\newcommand{\plus}{{+}}
\newcommand{\minus}{{-}}
\newcommand{\phiqp}{\phi_q^{\plus}}
\newcommand{\phiqm}{\phi_q^{\minus}}
\newcommand{\sss}{\mathcal S}
\newcommand{\pp}{\mathcal P}
\newcommand{\aaa}{\mathcal A}
\newcommand{\mm}{\mathcal M}
\newcommand{\cc}{\mathcal C}
\newcommand{\cl}{\operatorname{Cl}_2}
\newcommand{\ee}{{e}}
\newcommand{\xsup}{{S}}
\newcommand{\xinf}{{I}}

\newtheorem{theorem}{Theorem}
\newtheorem{lemma}{Lemma}
\newtheorem{proposition}{Proposition}
\newtheorem{corollary}{Corollary}
\newproclaim{definition}{Definition}
\newproclaim{remark}{Remark}

\makeatother

\begin{document}
\begin{frontmatter}

\title{On extrema of stable processes}
\runtitle{On extrema of stable processes}

\begin{aug}
\author[A]{\fnms{Alexey} \snm{Kuznetsov}\corref{}\thanksref{t1}\ead
[label=e1]{kuznetsov@mathstat.yorku.ca}}
\runauthor{A. Kuznetsov}
\affiliation{York University}
\address[A]{Department of Mathematics and Statistics\\
York University \\
Toronto, Ontario, M3J 1P3\\
Canada\\
\printead{e1}} 
\end{aug}

\thankstext{t1}{Supported by the Natural Sciences and Engineering Research
Council of Canada.}

\received{\smonth{1} \syear{2010}}
\revised{\smonth{6} \syear{2010}}

%
\begin{abstract}
We study the Wiener--Hopf factorization and the distribution of extrema
for general stable processes. By connecting the Wiener--Hopf factors
with a certain elliptic-like function we are able to obtain many
explicit and general results, such as infinite series representations
and asymptotic expansions for the density of supremum, explicit
expressions for the Wiener--Hopf factors and the Mellin transform of
the supremum, quasi-periodicity and functional identities for these
functions, finite product representations in some special cases and
identities in distribution satisfied by the supremum functional.
\end{abstract}

%
\setattribute{keyword}{AMS}{MSC2010 subject classification.}

\end{frontmatter}

\section{Introduction}\label{section_Introduction}

Stable processes are probably the most fascinating members of the vast
family of L\'evy processes. They enjoy the scaling property, which
states that the process $\{ X_{ct}\dvtx  t\ge0\}$ has the same
distribution as $\{c^{1/\alpha}X_t\dvtx  t\ge0\}$. Another
well-known process which also enjoys this property is Brownian motion,
for which we can compute explicitly the distribution (in many cases
even the joint distribution) of many important functionals; an
extensive collection of these facts can be found in \cite{Borodin}. In
the case of stable processes the situation is much less satisfactory:
we do not have general explicit results even for such a basic
functional as the supremum. However, there does exist an ``almost''
general result: Doney \cite{Doney1987} has obtained closed-form
expressions for the Wiener--Hopf factors for a dense set of parameters.
These results have inspired the main idea for our current work.

In \cite{Doney1987} it was proven that if parameters $\alpha$ and
$\rho
=\p(X_1>0)$
satisfy
%
%
\begin{equation}\label{eqn_def_Ckl}
\rho+k=\frac{l}{\alpha}
\end{equation}
for some integers $k$ and $l$, then there exists a fully explicit
representation for the Wiener--Hopf factors, given as a ratio of two
finite products of length $|k|$ and $|l|$. It is easy to see that when
$\alpha$ is irrational, relation (\ref{eqn_def_Ckl}) defines a
countable set of values of $\rho$, which is dense in $[0,1]$. Thus,
assuming that $\alpha$ is fixed and irrational, we have a function with
explicit representation for a dense countable set of values of the
parameter $\rho$. This situation reminds one of an analogous property
of elliptic functions: when a certain parameter is a rational number
there often exists an explicit expression, or at least a reduction to
some simpler functions. Therefore, given results obtained in
\cite{Doney1987}, it seems highly probable that there should exist a
connection between extrema of stable processes and elliptic functions.
The goal of this paper is to exhibit this connection and as a
consequence to obtain a remarkably large number of explicit results on
extrema of stable processes.

First we will present some definitions and notations. We will always
use the principal branch of the logarithm, which is defined in the
domain $|{\arg}(z)|<\pi$
by requiring that $\ln(1)=0$. Similarly, the power function will be
defined as $z^a=\exp(a\ln(z))$ in the domain $|{\arg}(z)|<\pi$.
We will work with a stable process $X_t$, which is a L\'evy process
with the characteristic exponent $\Psi(z)=-\ln(\e[\exp(\ii z X_1)])$
given by
%
%
\begin{equation}\label{def_Psi0}
\Psi(z)=c|z|^{\alpha} \biggl(1-\ii\beta\tan\biggl(\frac{\pi\alpha
}2\biggr) \operatorname{sgn}(z)\biggr),\qquad z\in\rr,
\end{equation}
where $c>0$, $\alpha\in(0,1)\cup(1,2)$ and $\beta\in[-1,1]$ (see
\cite{Bertoin} and \cite{Doney1987}).
The fact that the characteristic exponent is essentially a power
function of order $\alpha$
implies the scaling property, which states that
processes $\{X_{ct} \dvtx  t\ge0\}$ and $\{c^{1/\alpha}X_t \dvtx
t\ge0\}$ have the same distribution.

Following \cite{Doney1987} (and \cite{Zolotarev1957}) we rewrite
(\ref
{def_Psi0}) in a more convenient form. First of all, $c$ is just a
scaling parameter, thus without loss of generality we can assume that
$c^2=1+\beta^2 \tan(\frac{\pi\alpha}2)^2$.
Next we introduce a parameter
%
%
\begin{equation}\label{def_gamma}
\gamma=\frac{2}{\pi} \tan^{-1} \biggl(-\beta\tan\biggl(\frac{\pi\alpha
}2 \biggr)\biggr)
\end{equation}
and rewrite (\ref{def_Psi0}) as
%
%
\begin{equation}\label{def_Psi}
\Psi(z)=e^{\pi\ii\gamma/2} |z|^{\alpha} \mathbf{1}_{\{z>0\}
}+e^{-\pi\ii\gamma/2}|z|^{\alpha} \mathbf{1}_{\{z<0\}}.
\end{equation}

Another important parameter is $\rho=\p(X_1>0)$, which was computed in
closed form in \cite{Zolotarev1957}
%
%
\begin{equation}\label{def_rho}
\rho=\frac{1}2 \biggl(1-\frac{\gamma}{\alpha} \biggr)=\frac12+
\frac{1}{\pi\alpha} \tan^{-1} \biggl(\beta\tan\biggl(\frac{\pi\alpha}2
\biggr)\biggr).
\end{equation}
One can see that parameters $(\alpha,\beta)$ belong to the set
\[
\{\alpha\in(0,1), \beta\in(-1,1)\} \cup\{\alpha=1, \beta=0
\} \cup\{\alpha\in(1,2), \beta\in[-1,1]\},
\]
if and only if parameters $(\alpha,\rho)$ belong to the following set:
%
%
\begin{eqnarray}\label{def_set_A}
\aaa&=&\{\alpha\in(0,1), \rho\in(0,1)\} \cup\bigl\{\alpha=1, \rho
={ {\tfrac12}} \bigr\} \nonumber\\[-8pt]\\[-8pt]
&&{}\cup \{\alpha\in(1,2),
\rho\in[1-\alpha^{-1}, \alpha^{-1}]\}.\nonumber
\end{eqnarray}
We will call $\aaa$ the admissible set of parameters and will
parametrize the stable process $X_t$ by
$(\alpha,\rho)$. Note that we exclude the case when $\alpha\in(0,1)$
and $\rho=1$ $\{\rho=0\}$, as in this case the process $X_t$ $\{
-X_t\}$
is a subordinator and the Wiener--Hopf factorization is trivial. In the
limiting case as $\alpha\to2^-$ and $\rho\to1/2$ we obtain Brownian
motion $X_t=\sqrt{2} W_t$. When $\alpha\in(1,2)$ and $\rho=1-\alpha
^{-1}$ ($\rho=\alpha^{-1}$) the process $X_t$ is spectrally positive
(negative). In this case we have
complete information about the Wiener--Hopf factorization and
the distribution of extrema due to the work of
Bingham \cite{Bingham1975}, Doney \cite{Doney2008}, Bernyk, Dalang and
Peskir \cite{Bernyk2008} and Patie \cite{Patie2009}.

The following family of stable processes, which was introduced in \cite
{Doney1987}, will be very important in our paper.
\begin{definition}
For $k,l \in\zz$ define $\cc_{k,l}$ as the class of stable processes
$X_t$ with parameters $(\alpha,\rho) \in\aaa$ satisfying relation
(\ref
{eqn_def_Ckl}).
\end{definition}

Note that if $(\alpha,\rho)$ satisfy relation (\ref{eqn_def_Ckl}), then
$\alpha\in\q$ if and only if $\rho\in\q$.
When $\alpha\notin\q$ there exists a unique pair of integer numbers
$k,l$ such that relation (\ref{eqn_def_Ckl}) holds.
If $\alpha=\frac{m}{n}$ for some coprime integers $m,n$ and $X_t \in
\cc_{k,l}$,
then relation (\ref{eqn_def_Ckl}) holds for any pair $(\tilde k,
\tilde l)=(k+jn, l+jm)$, $j\in\zz$. In this case
we will assume that $0\le k < n$ and $1\le l < m$. Another important
fact is that the process $X_t$ with parameters $(\alpha,\rho)$ belongs
to $\cc_{k,l}$ [and $\alpha\in(1,2)$] if and only if the process
$\tilde X_t$ with parameters
$(\tilde\alpha, \tilde\rho)=(\alpha\rho,1/\alpha)$ belongs to
$\cc
_{-l,-k}$. This property, which can be easily verified using relation
(\ref{eqn_def_Ckl}), will appear in many formulas throughout this paper
and will be very useful for us later.
Another important property is that $X_t$ belongs to $\cc_{k,l}$ if and
only if the dual process $\hat X_t=-X_t$ belongs to $\cc_{-k-1,-l}$.
This can be easily checked using (\ref{eqn_def_Ckl}), as the dual
process has parameters $(\alpha,1-\rho)$.

Next, we define the supremum and infimum processes
\[
\xsup_t=\sup\{X_u \dvtx  0\le u \le t\},\qquad \xinf_t=\inf\{X_u \dvtx  0\le
u \le t\}.
\]
We introduce a random variable $\ee(q)\sim \operatorname{Exp}(q)$ (exponentially distributed
with parameter $q > 0$) which is independent of the process $X_t$.
We will use the following standard notation for the characteristic
functions of $\xsup_{\ee(q)}$ and
$\xinf_{\ee(q)}$, also known as the Wiener--Hopf factors:
\[
\phiqp(z)=\e[e^{\ii z \xsup_{\ee(q)}}],\qquad \phiqm(z)=\e
[e^{\ii z \xinf_{\ee(q)}}],\qquad z\in\rr.
\]
Note that since the random variable $\xsup_{\ee(q)}$ $\{\xinf_{\ee
(q)}\}$ is positive \{negative\}, function $\phiqp(z)$ $\{\phiqm(z)\}$
can be extended analytically into domain $\im(z)>0$ $\{\im(z)<0\}$.

To reduce the complexity of the problem, we note first that it is
enough to study the positive Wiener--Hopf factor $\phiqp(z)$, as the
corresponding information about $\phiqm(z)$ can be obtained by
considering the dual process. Second, the scaling property implies that
\[
\phiqp(z)=\phi^{\plus}_1(zq^{-1/\alpha});
\]
therefore we only need to consider the case $q=1$. This justifies our
choice of the following function as
the main object of study in this paper:
%
%
\begin{equation}\label{def_phi}
\phi(z)=\phi(z;\alpha,\rho)=\e[ e^{-z \xsup_{\ee(1)}} ],\qquad
\re(z) \ge0.
\end{equation}

This paper is organized as follows: in Section \ref{section_gen_props}
we establish a connection between $\phi(z)$
and a certain elliptic-like function $F(z;\tau)$, and we study various
analytical properties of the latter. In Section
\ref{section_rational_alpha} we derive an explicit expression (in terms
of the Clausen function) for the function $\phi(z)$ when $\alpha$ is rational.
In Section \ref{section_inf_product} we study the Wiener--Hopf
factorization in the general case:
we express the Wiener--Hopf factor $\phi(z)$ in terms of the
q-Pochhammer symbol for $\im(\alpha)>0$,
rederive some of the results obtained in \cite{Doney1987}, express
$\phi
(z)$ in terms of the
double gamma function and give a complete description of the
analytical properties of~$\phi(z)$.
In Section \ref{section_functional_eqns} we obtain a series representation
for $\ln(\phi(z))$ with very interesting convergence properties and
establish several functional equations satisfied by $\phi(z)$. In Section
\ref{section_mellin_transform} we study the Mellin transform of the
supremum $\xsup_1$: we obtain two
functional equations and establish quasi-periodicity of this function,
express it in terms of the double gamma function
and as a finite product when $X_t \in\cc_{k,l}$ and present two
distributional identities satisfied by
the supremum functional. Finally, in Section \ref{section_pdf_supremum}
we derive several convergent and asymptotic series representations for the
probability density function of the supremum $\xsup_1$.

\section{Connecting the Wiener--Hopf factors with elliptic
functions}\label{section_gen_props}

Our main tool in studying the Wiener--Hopf factor $\phi(z)$ will be a
certain function $F(z;\tau)$,
which has many properties similar to elliptic functions, but first we
need to define
two types of domains in the complex plane.
\begin{definition} Assume $\im(\tau)>0$. Define
\begin{eqnarray*}
\sss(\tau)&=&\{z \in\ccc\dvtx  |{\re}(z \bar\tau)|<\pi\im(\tau)\}, \\
\pp
(\tau)&=&\sss(\tau) \cap\{ z \in\ccc\dvtx  |{\im}(z)|<\pi\im(\tau)\}.
\end{eqnarray*}
\end{definition}

One can see that $\sss(\tau)$ is a strip in the complex plane which
contains $z=0$,
and $\pp(\tau)$ is a parallelogram if $\re(\tau)\ne0$.
In the case when $\re(\tau)=0$ we have $\sss(\tau)=\{z \dvtx  |{\im}(z)|<
\pi
\}$,
therefore $\sss(\tau)$ becomes a horizontal strip which does not depend
on $\tau$, while parallelogram $\pp(\tau)$
degenerates into the horizontal strip
$\{z \in\ccc\dvtx  |{\im}(z)|<\pi\min(\im(\tau),1) \}$. Domain $\pp
(\tau)$
satisfies the following important property: $z \in\pp(\tau)$ if and
only if
$(-\ii z/\tau) \in\pp(-1/\tau)$. Also note that $\sss(a \tau
)\equiv
\sss(\tau)$ for $a\in\rr^{\plus}$, therefore $\sss(\tau)$
depends only
on the argument of $\tau$.
\begin{definition} Assume $\im(\tau)>0$ and $z\in\sss(\tau)$. Define
%
%
\begin{equation}\label{def_F_z_tau}
F(z;\tau)=\int_{\rr} \frac{\dd x}{(1+e^{z+\ii\tau x})(1+e^{x})}.
\end{equation}
\end{definition}

The integrand in the definition (\ref{def_F_z_tau}) converges to zero
exponentially as $x\to+\infty$,
condition $\im(\tau)>0$ ensures that the same it true as $x\to
-\infty$
and condition $z\in\sss(\tau)$ guarantees that $\exp(z+\ii\tau
x)\ne
-1$ for $x \in\rr$.
Therefore the integral in (\ref{def_F_z_tau}) converges absolutely,
and the
function $F(z;\tau)$ is well defined. It is also clear that
this function is analytic in both variables $(z,\tau)\in\ccc^2$ as
long as $\im(\tau)>0$ and $z \in\sss(\tau)$.
\begin{proposition}\label{prop_main} Assume $(\alpha,\rho)\in\aaa$ and
$|{\operatorname{arg}} (z)|<\pi\min((1-\rho),\frac12)$. Then
%
%
\begin{eqnarray}\label{eq_lnphi_Fztau_main}
\frac{\dd}{\dd z} \ln(\phi(z))&=&\frac{1}{2\pi\ii z} \bigl[
F\bigl( \pi\ii\rho-\ln(z); \ii\alpha^{-1}\bigr) \nonumber\\[-8pt]\\[-8pt]
&&\hspace*{27.2pt}{}-
F\bigl( -\pi\ii\rho-\ln(z);\ii\alpha^{-1}\bigr) \bigr].\nonumber
\end{eqnarray}
\end{proposition}
\begin{pf}
We start with the following general integral representation for $\phiqp
(z)$, which was first derived in \cite{Mordecki} (Lemma 4.2) (another
proof was given in \cite{Kuznetsov}, Theorem 1(b)):
%
%
\begin{equation}\label{eq_lnphi_int0}
\ln(\phiqp(z))=\frac{z}{2\pi\ii} \int_{\rr} \ln\biggl(\frac
{q}{q+\Psi(u)}\biggr) \frac{\dd u}{u(u-z)}, \qquad \im(z)>0.
\end{equation}
As was shown in \cite{Mordecki}, this formula is valid for any L\'evy
process $X_t$, provided that the integral
\[
\int_{-\epsilon}^{\epsilon} \biggl|\frac{\Psi(u)}{u}\biggr|\,\dd u
\]
converges for some $\epsilon>0$, which is true in the case of stable
processes due to (\ref{def_Psi}).
Integral representation (\ref{eq_lnphi_int0}) is in fact equivalent to
Darling's integral (see \cite{Darling1956,Heyde1969}) which was
the main tool used in \cite{Doney1987}.

Next, we use our definition $\phi(z)=\phi^{\plus}_1(\ii z)$ and
formulas (\ref{def_Psi}) and (\ref{eq_lnphi_int0}) to find that
%
%
\begin{eqnarray}\label{eq_lnphi_int}
\ln(\phi(z))&=&-\frac{z}{2\pi} \biggl[\int_0^{\infty} \frac{\ln
(1+e^{\pi\ii\gamma/2}u^{\alpha})}{u(u-\ii z)}\,\dd u
\nonumber\\[-8pt]\\[-8pt]
&&\hspace*{27.7pt}{}+
\int_{0}^{\infty} \frac{\ln(1+e^{-\pi\ii\gamma/2}
u^{\alpha})}{u(u+\ii z)}\,\dd u\biggr],\qquad \re(z)>0.\nonumber
\end{eqnarray}
We assume first that $z \in\rr^+$ and obtain
%
%
\begin{eqnarray}\label{eq_lnphi_int2}
&&
\frac{\dd}{\dd z} \biggl[ z \int_0^{\infty} \frac{\ln(1+e^{
{\pi\ii\gamma}/2}u^{\alpha})}{u(u-\ii z)}\,\dd u \biggr]\nonumber\\
&&\qquad=
\int_0^{\infty} \frac{\ln(1+e^{\pi\ii\gamma/2}u^{\alpha
})}{(u- \ii z)^2}\,\dd u \nonumber\\[-8pt]\\[-8pt]
&&\qquad=- \int_0^{\infty} \ln(1+e^{{\pi\ii\gamma
}/2}u^{\alpha}) \,\frac{\dd}{\dd u} [ (u- \ii z)^{-1} ]
\,\dd u \nonumber\\
&&\qquad=
-\alpha\int_0^{\infty} \frac{u^{\alpha-1} }{(e^{-{\pi\ii
\gamma}/2}+u^{\alpha} )(u-\ii z)}\,\dd u,\nonumber
\end{eqnarray}
where in the last step we have applied integration by parts. Changing
the variable of integration $u=e^{-{y}/{\alpha}}$ in the last
integral in (\ref{eq_lnphi_int2}) we have
%
%
\begin{equation}\label{eq_lnphi_int3}
-\alpha\int_0^{\infty} \frac{u^{\alpha-1} }{(e^{-{\pi\ii
\gamma}/2}+u^{\alpha} )(u-\ii z)}\,\dd u=
\int_{\rr} \frac{\dd y}{(e^{-{y}/{\alpha}}-\ii z)(1+e^{-
{\pi\ii\gamma}/2+y})}.\hspace*{-32pt}
\end{equation}
The next step is to shift the contour of integration $\rr\mapsto\rr
+\frac{\pi\ii\gamma}2$ in the integral in the right-hand side of
(\ref
{eq_lnphi_int3}).
This step is justified, since the integrand decays exponentially
as $|{\re}(y)| \to\infty$, and the inequality $|\gamma|<\alpha$ [which
follows from (\ref{def_gamma})] guarantees that the integrand is
analytic in the horizontal strip $\im(y)<|\frac{\pi\gamma}2|$.
Thus, shifting the contour of integration and performing one
final change of variables $y=x+\frac{\pi i \gamma}2$ we obtain
%
%
\begin{eqnarray}\label{eq_lnphi_int4}
&&\int_{\rr} \frac{\dd y}{(e^{-y/\alpha}-\ii
z)(1+e^{-{\pi\ii\gamma}/2+y})}\nonumber\\
&&\qquad=
\int_{\rr+{\pi\ii\gamma}/2} \frac{\dd y}{(e^{-y/\alpha}-\ii
z)(1+e^{-{\pi\ii\gamma}/2+y})}\nonumber\\[-8pt]\\[-8pt]
&&\qquad =
\int_{\rr} \frac{\dd x}{(e^{-{x}/{\alpha}-{\pi\ii
\gamma}/({2\alpha})}-\ii z)(1+e^{x})}\nonumber\\
&&\qquad=
\frac{\ii}{z}\int_{\rr} \frac{\dd x}{(1+e^{{\pi\ii}/2
(1- {\gamma}/{\alpha})- \ln(z)-{x}/{\alpha}})(1+e^{x})}.\nonumber
\end{eqnarray}
Combining (\ref{eq_lnphi_int2}), (\ref{eq_lnphi_int3}), (\ref
{eq_lnphi_int4}) with the definitions of $F(z;\tau)$ and $\rho$ we
conclude that
\[
\frac{\dd}{\dd z} \biggl[ z \int_0^{\infty} \frac{\ln(1+e^{
{\pi\ii\gamma}/2}u^{\alpha})}{u(u-\ii z)}\,\dd u \biggr]
=\frac{\ii}{z} F\bigl( \pi\ii\rho-\ln(z);\ii\alpha^{-1} \bigr).
\]
In the case $z \in\rr^{\plus}$ equation (\ref{eq_lnphi_Fztau_main})
follows by taking the real part of both sides of the above equation and
using (\ref{eq_lnphi_int}). We can extend this result by analytic
continuation into domain $|{\arg} (z)|<\pi\min( (1-\rho
),\frac12)$,
since the right-hand side of (\ref{eq_lnphi_Fztau_main}) is analytic in
the region $|{\arg} (z)|<\pi(1-\rho)$ and due to (\ref
{eq_lnphi_int})
$\ln(\phi(z))$ is well defined and analytic in the domain
$|{\arg} (z)|<\frac{\pi}2$.
\end{pf}

At this point we would like to make several remarks which will be
important later. First of all,
the integral representation (\ref{eq_lnphi_int}) allows us to define
the function $\phi(z;\alpha,\rho)$ for all positive $\alpha$,
even though for $\alpha> 2$ these functions do not have an immediate
probabilistic interpretation.
Second, and most importantly, Proposition \ref{prop_main} allows us to
consider $\phi(z;\alpha,\rho)$ as a function which is analytic in all
three variables $(z,\alpha,\rho)$.
We do not need to specify explicitly the domain in $\ccc^3$ where
this function is analytic; for our purposes it is enough to know that for
$\epsilon>0$ small enough the function $\phi(z;\alpha,\rho)$ is
analytic in all three variables in the domain
$\{|z-1|<\epsilon, |{\arg}(\alpha)|<\epsilon, |\rho-\frac
12|<\epsilon\}
\subset\ccc^3$ [this follows\vspace*{1pt} immediately from (\ref
{eq_lnphi_Fztau_main}) and the domain of analyticity of $F(z;\tau)$].

In the next theorem we collect several properties of the function
$F(z;\tau)$.
As we will see later, Proposition \ref{prop_main} will allow us to
translate each of these properties
into an important statement about the Wiener--Hopf factor $\phi(z)$.
\begin{theorem}\label{thm_props_F_z_tau} Assume $\im(\tau)>0$.

\begin{enumerate}[(iii)]
\item[(i)] For $z\in\pp(\tau)$
%
%
\begin{equation}\label{F_1_over_tau}
F(z;\tau)=\frac{\ii}{\tau} F\biggl( \frac{\ii z}{\tau};-\frac{1}{\tau}
\biggr).
\end{equation}
\item[(ii)] For $z\in\sss(\tau)$ and $n\ge2$
%
%
\begin{equation}\label{F_ntau}
F(z; n \tau)=\frac{1}{n}\sum_{k=0}^{n-1} F\biggl( \frac1n\bigl(z+\pi
\ii(n-2k-1)\bigr); \tau\biggr).
\end{equation}
\item[(iii)] For $z \in\pp(\tau)$ and $0<\epsilon<\min(\frac{\pi
}2, \im
(-\frac{\pi}{2\tau}))$
%
%
\begin{equation}\label{F_as_integral2}
F(z;\tau)=\frac{1}{2} \int_{\rr+\ii\epsilon} \frac{e^{
{\ii z x}/{\pi}}}{\sinh(x)\sinh(\ii\tau x)} \,\dd x.
\end{equation}
\item[(iv)] For $z \in\sss(\tau) \cap-\sss(\tau)$
%
%
\begin{equation}\label{eq_F_minus_z}
F(z;\tau)=F(-z;\tau)-\frac{\ii z}{\tau}.
\end{equation}
\item[(v)] Assume $\re(\tau)\ne0$ and define $\delta=
\operatorname{sign}(\re(\tau))$. For $z\in\sss(\tau)$
%
%
\begin{eqnarray}\label{F_series1}
F(z;\tau)&=&\delta\sum_{k\ge0} \biggl[\frac{2\pi\ii}{1+\exp
(z+2\pi\delta(k+1/2) \tau)} \nonumber\\[-8pt]\\[-8pt]
&&\hspace*{25.5pt}{}+
\frac{2\pi\tau^{-1}}{1+\exp({\ii z}/{\tau}+{2\pi\delta
}/{\tau
}(k+1/2))} \biggr].\nonumber
\end{eqnarray}
\item[(vi)] Assume $\re(\tau)\ne0$. For $z\in\sss(\tau) \cap\{
z\dvtx  \re
(z)>0\}$
%
%
\begin{equation}\label{F_series2}
F(z;\tau)=-\pi\ii\sum_{k\ge1} (-1)^k
\biggl[\frac{e^{- k z}}{\sinh(\pi k \tau)} -\frac{\ii}{\tau}\frac{e^{-
{\ii k z}/{\tau}}}{\sinh(\pi k / \tau)} \biggr].
\end{equation}
\end{enumerate}
\end{theorem}
\begin{pf*}{Proof of Theorem \ref{thm_props_F_z_tau}}
(i) Assume $z>0$ and $\re(\tau)=0$. Then $\ii\tau^{-1} \in\rr^{\plus}$,
and performing
the change of variables $x=-\frac{\ii}{\tau} y +\frac{\ii z}{\tau
}$ in
the integral defining $F(z;\tau)$ [see (\ref{def_F_z_tau})] we
immediately obtain identity (\ref{F_1_over_tau}). The general case
can be obtained by analytic continuation, as both sides of
(\ref{F_1_over_tau}) are analytic in the domain $\{\im(\tau)>0, z\in
\pp(\tau)\}$.

(ii) We start with the following identity:
%
%
\begin{equation}\label{part_fractions_identity}
\frac{1}{1+w^n}=\frac{1}{n}\sum_{k=0}^{n-1} \frac{1}{1+we^{
{\pi\ii}/{n}(n-2k-1)}},
\end{equation}
which is just a partial fraction decomposition of the rational function
$(1+w^n)^{-1}$.
Again, let us assume that
$z\in\rr$ and $\re(\tau)=0$.
Applying this identity to (\ref{def_F_z_tau}) we obtain
\[
\int_{\rr} \frac{\dd x}{(1+e^{z+\ii n\tau x})(1+e^{x})}=
\frac{1}{n} \sum_{k=0}^{n-1} \int_{\rr} \frac{\dd
x}{(1+e^{{z}/n+{\pi\ii}/{n}(n-2k-1)+\ii\tau x})(1+e^{x})},
\]
which is equivalent to (\ref{F_ntau}). Using property $\sss(n\tau
)\equiv\sss(\tau)$ and the fact that
$z\in\sss(\tau)$ implies $\frac{1}{n}(z\pm\pi\ii(n-1))\in\sss
(\tau
)$, we find that both sides of (\ref{F_ntau})
are analytic in $\{\im(\tau)>0, z\in\sss(\tau)\}$, thus statement
(ii) follows by analytic continuation.

(iii) Assume $z \in\rr$ and $\re(\tau)=0$, the general
case can be established as usual by analytic continuation.
We rewrite (\ref{def_F_z_tau}) as
%
%
\begin{equation}\label{eq_F_z_tau_new}
F(z;\tau)=\frac{1}{4} e^{-{z}/2} \int_{\rr} \frac{e^{-(1+\ii
\tau){x}/2} }
{\cosh({x}/2 ) \cosh(({z+\ii\tau x})/2 )}
\,\dd x.
\end{equation}
Using formula 3.511.4 from \cite{Jeffrey2007} we find that for $y\in
\rr
$ and $|\tau|$ sufficiently small
\begin{eqnarray*}
\int_{\rr} \frac{e^{-(1+\ii\tau){x}/2+\ii x y} }
{\cosh({x}/2 ) } \,\dd x&=&-\frac{2\pi\ii}{\sinh(\pi
(y-{\tau}/2))}, \\
\int_{\rr} \frac{e^{-\ii x y}}
{ \cosh(({z+\ii\tau x})/2 )} \,\dd x&=&\frac{2\pi\ii}{\tau
} \frac{e^{{yz}/{\tau}}}{\cosh({\pi\ii y}/{\tau} )}.
\end{eqnarray*}
Applying Plancherel's theorem to (\ref{eq_F_z_tau_new}) and using the
above Fourier transform identities, we find that
\[
F(z;\tau)=\frac{\pi}{2\tau} \int_{\rr} \frac{e^{(y-
{\tau}/2){z}/{\tau} }}
{ \cosh({\pi\ii y}/{\tau} )\sinh(\pi(y-
{\tau}/2))}\,\dd y.
\]
To finish the proof, we shift the contour of integration in the above
integral $\rr\mapsto\rr+\frac{\tau}2-\ii\epsilon$ (where
$\epsilon>0$
is a sufficiently small number)
and perform the change of variables $y=\frac{\ii\tau}{\pi}x+\frac
{\tau}2$.

(iv) Let $f(x)$ be the integrand in (\ref{F_as_integral2}).
In order to derive identity (\ref{eq_F_minus_z}) we start with the
integral representation
(\ref{F_as_integral2}), shift the contour of integration $\rr+\ii
\epsilon\mapsto\rr- \ii\epsilon$ taking into account the residue at
$x=0$ and finally change the variable of integration $x=-y$
\begin{eqnarray*}
2F(z;\tau)&=&\int_{\rr+\ii\epsilon} f(x) \,\dd x=
-2\pi\ii \operatorname{Res}_{ x=0} f(x) +
\int_{\rr+\ii\epsilon} f(-y) \,\dd y \\ &=&-2\pi\ii
\operatorname{Res}_{ x=0} f(x)+2F(-z;\tau).
\end{eqnarray*}
The residue at zero of $f(x)$ is easily seen to be $\frac{z}{\pi\tau}$.

\mbox{\hphantom{i}}(v) Let us fix $n\ge1$ and assume $\re(\tau)>0$.
By shifting the contour of integration
in (\ref{def_F_z_tau}) $\rr\mapsto\rr-2\pi\ii n $ (while taking care
of the residues) and then
changing variable of integration $x=y-2\pi\ii n$, we find
\begin{eqnarray*}
F(z;\tau)&=&\sum_{k=0}^{n-1} \frac{2\pi\ii}{1+\exp(z+2\pi
(k+1/2) \tau)}\\
&&{}+
\sum_{k=0}^{n^*}
\frac{2\pi\tau^{-1}}{1+\exp({\ii z}/{\tau}+{2\pi}/{\tau
}(k+1/2))}
\\
&&{}+\int_{\rr} \frac{\dd y}{(1+e^{z+2n\pi\tau+\ii\tau y})(1+e^{y})},
\end{eqnarray*}
where $n^*=\max\{k\dvtx  \im(-\frac{\ii z}{\tau}-\frac{2\pi\delta
}{\tau
}(k+\frac12))<2\pi n\}$.
As $n\to+\infty$ we have $|e^{2n\pi\tau}| \to\infty$, which implies
that the integral in the above equation converges to zero, and since
$n^* \to\infty$ we obtain series expansion (\ref{F_series1}).
If $\re(\tau)<0$ the proof is identical, except that now we need to
shift the contour of integration
into the upper half-plane: $\rr\mapsto\rr+2\pi\ii n$.

(vi) The proof is essentially the same as in part (v):
we start with the
integral representation (\ref{F_as_integral2}) and shift the contour of
integration into the upper-half plane while taking care of the residues.
The details are left to the reader.
\end{pf*}

\section{Finite product representation when $\alpha$ is
rational}\label
{section_rational_alpha}

As a first illustration of the power of Proposition \ref{prop_main} and
Theorem \ref{thm_props_F_z_tau} we will derive
an explicit expression for $\phi(z)$ when $\alpha$ is rational and
$\rho
$ is completely arbitrary.
This expression involves the Clausen function, defined for $\theta\in
\rr$ as
%
%
\begin{equation}\label{def_cl2}
\cl(\theta)=\sum_{n \ge1} \frac{\sin(n \theta)}{n^2}.
\end{equation}
The Clausen function can also be defined as the imaginary part of
$\operatorname{Li}_2(e^{\ii\theta})$, where
$\operatorname{Li}_2(z)$ is the dilogarithm function (see \cite
{Lewin1981} for an extensive study of both dilogarithm
and Clausen functions).

Function $\cl(\theta)$ can be easily evaluated numerically. Definition
(\ref{def_cl2}) implies
that for $\theta\in\rr$ and $n \in{\mathbb Z}$
\[
\cl(\theta+2\pi n)=\cl(\theta),\qquad \cl(-\theta)=-\cl(\theta),\qquad
\cl(2\pi-\theta)=-\cl(\theta),
\]
and thus we only need to be able to compute this function
for $|\theta| < \pi$, where it can be done very efficiently with the
help of the following series
representation (see \cite{Borwein2000} and formula 4.28 in \cite{Lewin1981}):
\begin{eqnarray*}
\cl(\theta)&=&3\theta-\theta\ln\biggl(|\theta|\biggl(1-\frac{\theta
^2}{4\pi^2} \biggr)\biggr)
\\
&&{}-2\pi\ln\biggl(\frac{2\pi+\theta}{2\pi- \theta} \biggr)+\theta
\sum_{n\ge1} \frac{\zeta(2n)-1}{n(2n+1)} \biggl(\frac{\theta
}{2\pi} \biggr)^{2n},
\end{eqnarray*}
where $\zeta(s)$ is the Riemann zeta function. Note that $\zeta
(2n)-1=4^{-n}+9^{-n}+\cdots;$
thus when $|\theta|<\pi$ the terms in the above series decrease as
$O(n^{-2}16^{-n})$ and
we have a very fast convergence rate.
\begin{theorem}\label{thm_alpha_rational} Assume that $\alpha=\frac
{m}{n}$ where $m$ and $n$ are coprime natural numbers.
Define
\[
\theta=
\cases{
\cot^{-1} \bigl(\cot(\pi m \rho)+(-1)^{mn}z^m\sin(\pi m \rho)^{-1}
\bigr), &\quad if $m \rho\notin
\zz$,\cr
0, &\quad if $m \rho\in\zz$.}
\]
Then for $z>0$
%
%
\begin{eqnarray}\label{eq_phiqp_rational}\quad
\phi(z)&=& \exp\biggl(\frac{1}{2\pi mn} \bigl( \cl(2\theta)-\cl(2\pi m
\rho)-\cl(2\theta-2\pi m \rho)\bigr) \biggr) \nonumber\\
&&{}\times
\bigl(1+(-1)^{mn}2\cos(\pi m \rho) z^m +z^{2m} \bigr)^{-{\rho
}/({2n})}\nonumber\\[-8pt]\\[-8pt]
&&{}\times
\prod_{k=0}^{n-1}
\bigl(1+2\cos\bigl(\pi\alpha(\rho+2k+1) \bigr) z^{\alpha} +z^{2\alpha
}\bigr)^{({n-2k-1})/({2n})} \nonumber\\
&&{}\times
\prod_{j=0}^{m-1}
\biggl(1+2\cos\biggl({\frac{\pi}{\alpha}}(\alpha\rho+2j+1) \biggr)
z +z^{2}\biggr)^{({m-2j-1})/({2m})}.\nonumber
\end{eqnarray}
\end{theorem}

The proof of Theorem \ref{thm_alpha_rational} is based on the
following result:
\begin{lemma}\label{lemma_alpha_ratnl} If $m$ and $n$ are coprime
natural numbers and $|{\im}(z)|<\pi$
%
%
\begin{eqnarray}\label{eqn_Fztau_alpha_ratnl}
F\biggl(z;\ii\frac{m}{n}\biggr)&=&-\frac{n}{m} \frac
{z}{(-1)^{mn}e^{nz}+1} \nonumber\\
&&{}+ \sum_{k=0}^{n-1} \frac{{\pi\ii}/{n} (n-2k-1)}{\exp
(z+{\pi\ii m}/{n} (2k+1))+1} \\
&&{}+\frac{n}{m} \sum_{j=0}^{m-1} \frac{{\pi\ii}/{m}
(m-2j-1)}{\exp(z{n}/{m}+{\pi\ii n}/{m} (2j+1))+1}.\nonumber
\end{eqnarray}
\end{lemma}
\begin{pf}
The proof is in fact quite simple. The idea is to apply
(\ref{F_1_over_tau}) and (\ref{F_ntau}) in order to transform
$F(z;\ii
\frac{m}{n} )$ into a sum
of $F(\cdot;\ii)$, which can be evaluated explicitly. First we use
(\ref
{F_ntau}) and obtain
\[
F\biggl(z;\ii\frac{m}n \biggr)= \frac{1}{m} \sum_{j=0}^{m-1}
F\biggl( \frac{1}{m}\bigl(z+\pi\ii(m-2j-1)\bigr);\frac{\ii}{n} \biggr).
\]
Next, applying (\ref{F_1_over_tau}) to each function $F(\cdot, \frac
{\ii
}{n})$ in the right-hand side of the above formula we find that
\[
F\biggl(z;\ii\frac{m}n \biggr)= \frac{n}{m} \sum_{j=0}^{m-1}
F\biggl( \frac{n}{m}\bigl(z+\pi\ii(m-2j-1)\bigr);\ii n \biggr).
\]
Again, we apply (\ref{F_ntau}) to each function $F(\cdot, \ii n)$ in
the right-hand side of the above formula and deduce
%
%
\begin{equation}\label{lemma_artn_proof1}
F\biggl(z;\ii\frac{m}{n}\biggr)=\frac{1}{m} \sum_{k=0}^{n-1} \sum
_{j=0}^{m-1}
F\biggl(\frac{z}{m}+\frac{\pi\ii}n (n-2k-1)+\frac{\pi\ii}{m}
(m-2j-1);\ii\biggr).\hspace*{-32pt}
\end{equation}
Now we only need to evaluate $F(z;\ii)$, which can be done as follows:
\begin{eqnarray*}
F(z;\ii)&=&\frac{e^{-{z}/2}}4\int_{\rr} \frac{\dd x}{\cosh
(({x-z})/2) \cosh({x}/{2})}\\
&=&
\frac{e^{-{z}/2}}2\int_{\rr} \frac{\dd x}{\cosh(x-
{z}/2)+ \cosh({z/2})}\\
&=&
e^{-{z}/2}\int_0^{\infty} \frac{\dd x}{\cosh(x)+
\cosh({z/2})}=
\frac{z}{e^z-1},
\end{eqnarray*}
where in the last step we have used formula 3.514.1 from \cite
{Jeffrey2007}. Finally, the double sum in the right-hand side of
(\ref{lemma_artn_proof1}) can be reduced to
single sums in (\ref{eqn_Fztau_alpha_ratnl}) using the following identity:
\[
\frac{1}{1+w^n}=\frac{1}{n}\sum_{k=0}^{n-1} \frac{1}{1+we^{
{\pi\ii m}/{n}(n-2k-1)}},
\]
which is just a different way of writing the partial fraction
decomposition (\ref{part_fractions_identity}).
\end{pf}
\begin{pf*}{Proof of Theorem \ref{thm_alpha_rational}}
Using Proposition \ref{prop_main} and formula (\ref
{eqn_Fztau_alpha_ratnl}) we find that
\begin{eqnarray*}
\frac{\dd}{\dd z} \ln(\phi(z))&=&(-1)^{mn+1}\frac{\sin(\pi m \rho
)}{\pi
m n}
\frac{\ln(z^m)mz^{m-1}}{1+(-1)^{mn}2\cos(\pi m \rho) z^m
+z^{2m}}\\
&&{}-
\frac{\rho}{2n} \frac{(-1)^{mn}2\cos(\pi m \rho) m z^{m-1}
+2mz^{2m-1}}{1+(-1)^{mn}2\cos(\pi m \rho) z^m +z^{2m}}\\
&&{}+\sum_{k=0}^{n-1} \biggl(\frac{n-2k-1}{2n} \biggr)
\frac{2\cos(\pi\alpha(\rho+2k+1) ) \alpha z^{\alpha-1}
+2\alpha z^{2\alpha-1}}{1+2\cos(\pi\alpha(\rho+2k+1) )
z^{\alpha} +z^{2\alpha}}\\
&&{}+\sum_{j=0}^{m-1} \biggl(\frac{m-2j-1}{2m} \biggr)
\frac{2\cos({\pi}/{\alpha}(\alpha\rho+2j+1) )
+2z}{1+2\cos({\pi}/{\alpha}(\alpha\rho+2j+1) ) z +z^{2}}.
\end{eqnarray*}
To complete the proof we just need to integrate the above identity and
use the fact that $\phi(0)=1$. Most of the integrals are obvious,
except for the following one:
\[
\delta\sin(\beta)\int_0^{z}\frac{\ln(u)}{1+2\delta\cos(\beta)
u+u^2} \,\dd u=
\frac{1}2 [ \cl(2\theta-2\beta)-\cl(2\theta)+\cl(2\beta) ],
\]
where $\delta\in\{-1, 1\}$ and $
\theta=\cot^{-1} ( \cot(\beta)+\delta z \sin(\beta)^{-1} )$.
However, this integral is just a particular case
of formula 8.18 in \cite{Lewin1981}, and we refer the interested reader
to this great book for all the details.
\end{pf*}

Note that while Theorem \ref{thm_alpha_rational} is different from the
result obtained in \cite{Doney1987}, there
are also some similarities. In particular, both results give explicit
formulas for a dense uncountable set of parameters
$(\alpha,\rho)$. If $\rho$ is also a rational number of the form
$\rho
=\frac{k}m$, then $\theta=0$, and formula (\ref{eq_phiqp_rational})
does not contain function $\cl$.
It is easy to see that the set of points
\[
\biggl\{\alpha=\frac{m}{n}, \rho=\frac{k}m \biggr\}
\]
is dense in the set of all admissible parameters $\aaa$. It is probably
true that in this case formula (\ref{eq_phiqp_rational})
can be reduced to expression (1.11) in \cite{Doney1987}:
due to\vspace*{1pt} the fact that the numbers $m$ and $n$ are coprime we can find integers
$\tilde k$ and $\tilde l$, such that $k=\tilde l n-\tilde k m$, thus
$\rho+\tilde k = \tilde l / \alpha$ and the process
$X_t$ belongs to the class $\cc_{\tilde k, \tilde l}$.

\section{Infinite product representation and analytic
continuation}\label{section_inf_product}
In this section we will study analytic properties of the Wiener--Hopf
factor $\phi(z)$. The most important results are the
representations in terms of the q-Pochhammer symbol and the double
gamma function.
\begin{definition}\label{def_q_Pochhammer} For $n \in{\mathbb N}$ we
define the q-Pochhammer symbol as
\[
(a;q)_n=\prod_{k=0}^{n-1} (1-aq^k)
\]
and $(a;q)_0=1$. If $|q|<1$ we define $(a;q)_{\infty}=\prod
_{k\ge0} (1-aq^k)$.
\end{definition}
\begin{theorem}\label{thm_phi_q_product}
Assume that $\im(\alpha)>0$ and define
%
%
\begin{equation}\label{def_q_tilde_q}
q=e^{2\pi\ii\alpha}, \qquad \tilde q=e^{-{2\pi\ii}/{\alpha}}.
\end{equation}
Then for $|z|<\min\{\sqrt{|q|},\sqrt{|\tilde q|}\}$ and $|{\arg}
(z)|<\pi$
%
%
\begin{equation}\label{eq_phi_q_product}
\phi(z)=\frac{(- z \sqrt{\tilde q} e^{-\pi\ii\rho};\tilde q
)_{\infty}
(-z^{\alpha} \sqrt{q} e^{\pi\ii\rho\alpha}; q)_{\infty}}
{(- z \sqrt{\tilde q} e^{\pi\ii\rho};\tilde q)_{\infty}
(-z^{\alpha} \sqrt{q} e^{-\pi\ii\rho\alpha}; q)_{\infty}}.
\end{equation}
\end{theorem}
\begin{pf}
Using Proposition \ref{prop_main} and series expansion (\ref
{F_series1}) we find
\begin{eqnarray*}
\frac{\dd}{\dd z} \ln(\phi(z))&=&\sum_{k\ge0} \biggl[
\frac{1}{z+e^{\pi\ii\rho} \tilde q^{-k-1/2}}-\frac{\alpha
z^{\alpha-1}}{z^{\alpha}+e^{\pi\ii\rho\alpha} q^{-k-1/2}}
\\
&&\hspace*{18.2pt}{}-\frac{1}{z+e^{-\pi\ii\rho} \tilde q^{-k-1/2}}+\frac
{\alpha
z^{\alpha-1}}{z^{\alpha}+e^{-\pi\ii\rho\alpha} q^{-k-1/2}}
\biggr].
\end{eqnarray*}
To complete the proof we only need to integrate both sides of the above
identity and use the fact that $\phi(0)=1$.
\end{pf}

As a corollary of Theorem \ref{thm_phi_q_product} we can derive some
results obtained in \cite{Doney1987}.
\begin{corollary}
If $X_t \in{\mathcal C}_{k,l}$, then for $|{\arg}(z)|<\pi$
\[
\phi(z)=
\cases{
\dfrac{(z^{\alpha} (-1)^{1-l} q^{({1-k})/2};q)_k}
{(z (-1)^{1-k} \tilde q^{({1-l})/2};\tilde q)_l},
&\quad if
$l>0$,\vspace*{2pt}\cr
\dfrac{(z (-1)^{1+k} \tilde q^{({1+l})/2};\tilde q)_{|l|}}
{(z^{\alpha} (-1)^{1+l} q^{({1+k})/2};q)_{|k|}},
&\quad if
$l<0$.}
\]
\end{corollary}
\begin{pf}
Assume first that $\im(\alpha)>0$. We use (\ref{eq_phi_q_product}) and
the following identity for the q-Pochhammer symbol:
\[
\frac{(a;q)_{\infty}}{(aq^n;q)_{\infty}}=(a;q)_n
\]
(which can be easily obtained from Definition \ref{def_q_Pochhammer})
to derive the above finite product representation. The last step is
to let $\im(\alpha) \to0^+$ and use analytic continuation.
\end{pf}

A major disadvantage of Theorem \ref{thm_phi_q_product} is that we can
not extend it
to the case that interests us, when $\alpha$ is a real number, since
this would imply that $|\tilde q|=|q|=1$
and the infinite products do not converge. Therefore our next goal is
to transform
(\ref{eq_phi_q_product}) into an expression which is well
defined for both real and complex values of $\alpha$, and as we will
see later,
the double gamma function will be the magic key on our quest towards
analytic continuation.

Let us give some intuition on how one might guess that the double gamma
function is in fact the right tool to use.
Consider the following function:
\[
f(w)=\frac{(- e^{w} \sqrt{\tilde q} e^{-\pi\ii\rho};\tilde q
)_{\infty}}{(-e^{\alpha w} \sqrt{q} e^{-\pi\ii\rho\alpha}; q
)_{\infty}},
\]
which enters into formula (\ref{eq_phi_q_product}) with a change of
variables $z=e^{w}$. Using Definition~\ref{def_q_Pochhammer} of the
q-Pochhammer symbol we find that the numerator in the above expression
for $f(w)$ has simple zeros at points $w=w_{m,n}+\pi\ii\rho$ for
$m\ge0$ and $n \in\zz$, while the denominator has simple zeros at
$w=w_{m,n}+\pi\ii\rho$ for $m\in\zz$ and $n < 0$, where
%
%
\begin{equation}\label{def_w_mn}
w_{m,n}=\pi\ii[(2m+1)\alpha^{-1}+(2n+1) ],\qquad m,n\in\zz.
\end{equation}
We see that ``half'' of the zeros of the numerator coincide with
``half'' of the zeros of the denominator, and thus $f(w)$ has simple
zeros in the first quadrant of the lattice $w_{m,n}+\pi\ii\rho$ (where
$m\ge0$ and $n\ge0$) and simple poles in the third quadrant
$w_{m,n}+\pi\ii\rho$ ($m<0$ and $n<0$). Therefore it might be possible
to separate the function $f(w)$ into the ratio of two entire functions,
if we could only find a function which has simple zeros in the given
quadrant of the lattice. As it turns out, such a function was
introduced by Alexeiewsky in 1889 and, luckily enough, all of its
properties which are important for us were established by Barnes in
\cite{Barnes1899} (see also~\cite{Barnes1901}).

We will follow the original paper \cite{Barnes1899},
where Barnes introduces the double gamma function as an infinite
product in Weierstrass's form: for all
$|{\arg}(\tau)|<\pi$ and $z\in\ccc$
%
%
\begin{eqnarray}\label{def_G_z_tau}
G(z;\tau)&=&\frac{z}{\tau} e^{a{z}/{\tau}+b{z^2}/({2\tau})}
\nonumber\\[-8pt]\\[-8pt]
&&{}\times\prod
_{m\ge0} {\prod_{n\ge0}} ^{\prime}
\biggl(1+\frac{z}{m\tau+n} \biggr)e^{-{z}/({m\tau+n})+
{z^2}/({2(m\tau+n)^2})}.\nonumber
\end{eqnarray}
Here the prime in the second product means that the term corresponding
to $m=n=0$ is omitted.
Note that by definition $G(z;\tau)$ is an entire function in $z$ and
has simple zeros on the lattice $m\tau+n$, $m\le0$, $n\le0$.
Barnes proves that $G(z;\tau)$ can also be expressed as a single
infinite product of gamma functions
%
%
\begin{equation}\label{eq_G2_inf_prod_Gamma}\quad
G(z;\tau)=\frac{1}{\tau\Gamma(z)} e^{\tilde a {z}/{\tau
}+\tilde
b{z^2}/({2\tau^2})}
\prod_{m\ge1} \frac{\Gamma(m\tau)}{\Gamma(z+m\tau)} e^{z\psi
(m\tau)+{z^2}/2 \psi'(m\tau)},
\end{equation}
where $\psi(z)=\frac{\dd}{\dd z}\ln(\Gamma(z))$ is the digamma function
(see \cite{Jeffrey2007}). The constants $\tilde a$ and $ \tilde b$ are
related to $a$ and $b$ as follows:
\[
\tilde a=a-\gamma\tau,\qquad \tilde b=b+\frac{\pi^2 \tau^2}6,
\]
where $\gamma=-\psi(1)$ is the Euler--Mascheroni constant. One of the
most important properties of the double gamma function is
that it is quasi-periodic with periods $1$ and $\tau$, that is,
%
%
\begin{eqnarray}\label{funct_rel_G}
G(z+1;\tau)&=&\Gamma\biggl(\frac{z}{\tau}\biggr) G(z;\tau), \nonumber\\[-8pt]\\[-8pt]
G(z+\tau;\tau)&=&(2\pi)^{({\tau-1})/2}\tau^{-z+1/2}
\Gamma(z) G(z;\tau),\nonumber
\end{eqnarray}
provided that constants $\tilde a$ and $\tilde b$ are chosen in a
particular way
%
%
\begin{eqnarray}\label{def_a_b_tilde}
\tilde a&=&\frac{\tau}2 \ln(2\pi\tau)+\frac12 \ln(\tau)-\tau
C(\tau),\\
\tilde b&=&-\tau\ln(\tau)-\tau^2 D(\tau).
\end{eqnarray}
Here $C(\tau)$ and $D(\tau)$ are certain transcendental functions of
$\tau$ which can be computed as the following limits as $m\to+\infty$:
\begin{eqnarray*}
C(\tau)&=&\sum_{k=1}^{m-1}\psi(k\tau)+\frac12 \psi(m\tau)-\frac
{1}{\tau}\ln\biggl(\frac{\Gamma(m\tau)}{\sqrt{2\pi}}\biggr)
\\ &&{}-\frac{\tau}{12}\psi'(m\tau)+\frac{\tau^3}{720}\psi
^{(3)}(m\tau
)+O(m^{-5}), \\
D(\tau)&=&\sum_{k=1}^{m-1}\psi'(k\tau)+\frac12 \psi'(m\tau
)-\frac
{1}{\tau}\psi(m\tau)\\ &&{}-
\frac{\tau}{12}\psi''(m\tau)+\frac{\tau^3}{720}\psi^{(4)}(m\tau
)+O(m^{-6}).
\end{eqnarray*}
It turns out that with this choice of constants we also have $G(1;\tau
)=1$ (see~\cite{Barnes1899}).
There exists a different and slightly simpler expression for these
constants (see~\cite{Lawrie1994}),
however we have decided to use the original Barnes formulas as they are
more convenient for numerical calculations.
It is also possible to give an integral representation for $\ln
(G(z;\tau
))$ (see \cite{Lawrie1994}) and several asymptotic expansions (see
\cite
{Bill1997}).

The following two facts about the double gamma function will be very
important for us later. The first result was derived in \cite{Barnes1899}
and it is an analogue of the reflection formula for the gamma function
%
%
\begin{equation}\label{eq_G_tau_q_Poch}
-2\pi\ii\tau G\biggl(\frac12+z;\tau\biggr)G\biggl(\frac
{1}{2}-z;-\tau\biggr)
=\frac{(-e^{2\pi\ii z};q)_{\infty}}{(q;q)_{\infty}},
\end{equation}
where $q=e^{2\pi\ii\tau}$ and $\im(\tau)>0$.
The second result is the following transformation formula:
%
%
\begin{equation}\label{eq_G_1_over_tau}
G(z;\tau)=(2\pi)^{{z}/2 (1-1/{\tau} )} \tau^{-
{z^2}/({2\tau})+{z/2} ( 1+ {1/\tau} )-1 } G
\biggl(\frac{z}{\tau};\frac{1}{\tau}\biggr).
\end{equation}
We were not able to find any reference for this result; however, it is
not hard to prove it directly.
It follows from the definition (\ref{def_G_z_tau}) that the zeros of
$G(z;\tau)$ coincide with the zeros of the function on
the right-hand side of (\ref{eq_G_1_over_tau}). Performing some
straightforward but tedious computations one can check
that the latter function satisfies both functional equations (\ref
{funct_rel_G}).
Thus their ratio must be a function with two periods $1$ and $\tau$,
which is analytic in the entire complex plane, therefore
assuming that $\tau\notin\rr$ we conclude that this ratio must be a
constant. The fact
that the value of the constant is one can be established with the help
of (\ref{funct_rel_G}). We have included the detailed
proof of formula (\ref{eq_G_1_over_tau}) in the Appendix of the online
version of this paper, see \href
{http://arxiv.org/abs/1001.0991}{arXiv:1001.0991}.

Now we are ready to state and prove our main result in this section.
\begin{theorem}\label{thm_phi_Gamma2}
For $\alpha>0$, $\rho\in(0,1)$ and $|{\arg}(z)|<\pi$
%
%
\begin{eqnarray}\label{eq_phi_Gamma2}
\phi(z)&=&\bigl( 2 \pi\sqrt{z} \bigr)^{-\alpha\rho}
\frac{G(1/2 +{\alpha}/2 (1+\rho+{\ln(z)}/({\pi\ii
})) ; \alpha) }
{G( 1/2 +{\alpha}/2 (1-\rho+{\ln(z)}/({\pi\ii
})) ; \alpha) } \nonumber\\[-8pt]\\[-8pt]
&&{}\times
\frac{G(1/2 +{\alpha}/2 (1+\rho-{\ln(z)}/({\pi\ii
})) ; \alpha) }
{G( 1/2 +{\alpha}/2 (1-\rho-{\ln(z)}/({\pi\ii
})) ; \alpha) }.\nonumber
\end{eqnarray}
\end{theorem}
\begin{pf}
Assume that $\im(\alpha)>0$. We start with formula (\ref{eq_phi_q_product})
and apply identity (\ref{eq_G_tau_q_Poch}) to each q-Pochhammer symbol
to obtain
%
%
\begin{eqnarray}\label{thm4_proof0}
\phi(z)&=&\frac{G(1/2 -{1}/({2\alpha})-{\rho}/2+
{\ln(z)}/({2\pi\ii}) ; -{1}/{\alpha} ) }
{G( 1/2 -{1}/({2\alpha})+{\rho}/2+{\ln(z)}/({2\pi
\ii
}) ; -{1}/{\alpha} ) }\nonumber\\
&&{}\times\frac{G(1/2 +{1}/({2\alpha})+{\rho}/2-{\ln
(z)}/({2\pi
\ii}) ; {1}/{\alpha} ) }
{G( 1/2 +{1}/({2\alpha})-{\rho}/2-{\ln(z)}/({2\pi
\ii
}) ; {1}/{\alpha} ) }\nonumber\\[-8pt]\\[-8pt]
&&{}\times
\frac{G(1/2 + {\alpha}/{2}+{\rho\alpha}/2+\alpha
{\ln(z)}/({2\pi\ii}) ; \alpha) }
{G( 1/2 + {\alpha}/{2}-{\rho\alpha}/2+\alpha
{\ln
(z)}/({2\pi\ii}) ; \alpha) }\nonumber\\
&&{}\times
\frac{G(1/2 -{\alpha}/{2}-{\rho\alpha}/2-\alpha
{\ln(z)}/({2\pi\ii}) ; -\alpha) }
{G( 1/2 - {\alpha}/{2}+{\rho\alpha}/2-\alpha
{\ln
(z)}/({2\pi\ii}) ; -\alpha) }.\nonumber
\end{eqnarray}
Next, we transform the four double gamma functions in the first ratio
(the ones having $\tau=\pm\alpha^{-1}$) using identity
(\ref{eq_G_1_over_tau}) and after simplifying the resulting formula, we
obtain (\ref{eq_phi_Gamma2}). These calculations
are straightforward but lengthy and tedious; the interested reader can
find them in the Appendix of the online version of this paper, see
\href
{http://arxiv.org/abs/1001.0991}{arXiv:1001.0991}.

Now that we have established (\ref{eq_phi_Gamma2}) for the case $\im
(\alpha)>0$, we note that $G(z;\tau)$ is analytic and well defined when
$\tau\in\rr^{\plus}$.
Thus by analytic continuation
(\ref{eq_phi_Gamma2}) gives us an expression which is valid for
$\alpha
\in\rr^{\plus}$.
\end{pf}

Now we can give a complete description of the analytic structure of
$\phi(z)$.
\begin{corollary}\label{corollary_phi_analytic}
Assume that $\alpha>0$ and $\rho\in(0,1)$. The function $\phi(z)$ can
be analytically continued into the domain $\{z \in\ccc\dvtx  |{\arg}
(z)|<\pi\}
$, except that it has simple poles at points $\{-\exp(\pm\pi\ii
(\rho
-\alpha^{-1})\}$ when $\alpha\rho>1$. The function $\phi(z)$ has a
branch point at $z=0$ while the function
$w \mapsto\phi(e^w)$ is meromorphic in $\ccc$ and quasiperiodic with
periods $2\pi\ii$ and $2\pi\ii\alpha^{-1}$
%
%
\begin{eqnarray}\label{phi_analytic_ctnn}
\phi(e^{w+2\pi\ii})&=&\phi(e^w) \frac
{1+e^{\alpha w +\pi\ii\alpha(1-\rho)}}{1+e^{\alpha w + \pi\ii
\alpha
(1+\rho)}},
\\
\phi(e^{w+{2\pi\ii}/{\alpha}})&=&\phi(e^w)
\frac{1+ e^{w+\pi\ii(\alpha^{-1}-\rho)}}{1+ e^{w+ \pi\ii(\alpha
^{-1}+\rho)}}.\nonumber
\end{eqnarray}
The function $\phi(e^w)$ has roots at
\[
\{w_{m,n}+\pi\ii\rho\}_{m\ge0, n\ge0} \cup\{w_{m,n}-\pi\ii\rho
\}
_{m < 0, n<0}
\]
and poles at
\[
\{w_{m,n}-\pi\ii\rho\}_{m \ge0, n \ge0} \cup\{w_{m,n}+\pi\ii
\rho
\}_{m< 0, n<0},
\]
where $w_{m,n}=\pi\ii[(2m+1)\alpha^{-1}+(2n+1) ]$. All
roots and poles are simple if $\alpha$ is irrational.
\end{corollary}
\begin{pf}
The analytic continuation result follows from expression (\ref
{eq_phi_Gamma2}), as the double gamma function $G(z;\tau)$ has simple
zeros at points $z=-(m\tau+n)$, $n,m \ge0$. The quasi-periodicity
expressions and formulas for roots/poles of $\phi(e^w)$
follow at once from the infinite product representation (\ref
{eq_phi_q_product}).
\end{pf}

Note that if $(\alpha,\rho)\in\aaa$, then $\alpha\rho\le1$, which
implies that $\phi(z)$ is analytic in $|{\arg}(z)|<\pi$.
This result coincides with analytic continuation result for $\phi(z)$
when $X_t \in\cc_{k,l}$, which was obtained in \cite{Doney1987}.
Also, Corollary \ref{corollary_phi_analytic} gives us an insight into
the mysterious relation (\ref{eqn_def_Ckl}), which
was the key to obtain explicit expressions for $\phi(z)$ in \cite
{Doney1987}. Assuming $\rho=l \alpha^{-1} -k$ we find that all zeros
and poles of $\phi(e^w)$ lie on the same lattice $w_{m,n}$, and most of
them are canceled allowing us to express $\phi(e^w)$ as a ratio of
finite products of
trigonometric functions.

Applying the infinite product representation (\ref
{eq_G2_inf_prod_Gamma}) to Theorem \ref{thm_phi_Gamma2},
we obtain the following expression for $\phi(z)$ as an infinite product
of gamma functions. The detailed proof can be found in the Appendix of
the online version of this paper, see \href
{http://arxiv.org/abs/1001.0991}{arXiv:1001.0991}.
\begin{corollary}\label{corollary_phi_inf_prod_gamma} For $\alpha>0$,
$\rho\in(0,1)$ and $|{\arg}(z)|<\pi$
\begin{eqnarray*}
\phi(z)&=&z^{-{\alpha\rho}/2} e^{-\alpha\rho(2C(\alpha
)+(\alpha
+1)D(\alpha))}\\
&&{}\times\prod_{m\ge0}
e^{\alpha\rho(2\psi(m\alpha)+(\alpha+1) \psi'(m\alpha))\mathbf{1}\{
m >
0\}}\\
&&\hspace*{28pt}{}\times
\biggl(\Gamma\biggl( \frac12 +\frac{\alpha}2 \biggl(2m+1-\rho+\frac{\ln
(z)}{\pi\ii}\biggr) \biggr)\\
&&\hspace*{28pt}{}\times\Gamma\biggl( \frac12 +\frac{\alpha}2 \biggl(2m+1-\rho-\frac{\ln(z)}{\pi
\ii}\biggr)\biggr)\biggr)\\
&&\hspace*{28pt}{}\times\biggl(\Gamma\biggl( \frac12 +\frac{\alpha}2 \biggl(2m+1+\rho+\frac{\ln
(z)}{\pi\ii}\biggr) \biggr)\\
&&\hspace*{28pt}{}\times\Gamma\biggl( \frac12 +\frac{\alpha}2 \biggl(2m+1+\rho-\frac{\ln(z)}{\pi
\ii}\biggr) \biggr)\biggr)^{-1}.
\end{eqnarray*}
\end{corollary}

\section{Series representation and functional equations}\label
{section_functional_eqns}

As we have seen in the previous sections, rationality of certain
parameters plays an important role in our results.
In our next result there appear even more intriguing connections with
Number Theory.
It turns out that our next result is valid for almost all parameters
$\alpha$,
except for rational numbers and for those real numbers which can be
approximated by rational numbers unusually well.
\begin{definition}
A real number $x \in\rr\setminus\q$ belongs to the set ${\mathcal
L}$ if there exist $b>1$ and $c>0$ such that the inequality
\[
\biggl| \alpha-\frac{p}{q} \biggr| < \frac{c}{b^{q}}
\]
is satisfied for infinitely many coprime integers $p$ and $q$.
\end{definition}

The set ${\mathcal L}$ is a proper subset of the set of Lioville
numbers, which satisfy the following, weaker property: for all $\mu>0$
there exists $c>0$
such that the inequality
\[
\biggl| x-\frac{p}{q} \biggr| < \frac{c}{q^{\mu}}
\]
is satisfied for infinitely many coprime integers $p$ and $q$. A
celebrated result by Liouville states
that any algebraic number is not a Liouville number, but this is also
true for many other numbers.
In fact, almost every number is not a Liouville number,
as the set of Liouville numbers, while being dense in $\rr$, has zero
Lebesgue measure (see Theorem 32 in \cite{Khinchin}). The same is true
for the set ${\mathcal L}$: it is closed under multiplication by
rational numbers, and therefore it is dense in $\rr$, yet it has zero
Lebesgue measure.
\begin{theorem}\label{thm_lnphi_series}
Let $\alpha>0$ and $\rho\in(0,1)$. Assume that $\alpha\notin
{\mathcal L}\cup\q$.
Then if $|z|<1$ and $|{\arg}(z)| < \pi$
%
%
\begin{eqnarray}\label{eq_lnphi_series}
\phi(z)&=&\exp\biggl[ \sum_{k\ge1} \frac{\sin(\pi k \rho)}
{k\sin({\pi k }/{\alpha})}
(-1)^k z^k \nonumber\\[-8pt]\\[-8pt]
&&\hspace*{19.42pt}{}+ \sum_{k\ge1} \frac{\sin(\pi k \alpha
\rho)}{k\sin(\pi k \alpha)}
(-1)^k z^{\alpha k} \biggr].\nonumber
\end{eqnarray}
\end{theorem}
\begin{pf} Assume first that $z$ is a real number, such that $|z|<1$
and \mbox{$\im(\alpha)>0$}. Then $\re(\pm\pi\ii\rho- \ln(z))>0$ and
formula (\ref{eq_lnphi_series}) can be easily obtained using
Proposition \ref{prop_main} and
series expansion (\ref{F_series2}).
The hard part is to extend validity of (\ref{eq_lnphi_series}) to the
case when $\alpha$ is real,
and this is not a trivial matter. For example, when $\alpha$ is a
rational number the series on the right-hand side of
(\ref{eq_lnphi_series}) is not defined, as some terms will include
division by zero.

Assume that $\im(\alpha)>0$ and $\re(\alpha)$ is positive, irrational
and does not belong to~${\mathcal L}$. Then using the above definition
of the
set ${\mathcal L}$, we find that for every $b>1$ there exists $c>0$
such that $|k \re(\alpha)-n|>c k b^{-k}$ for all integers $k$ and
$n$. Therefore
we have
\[
|{\sin}(\pi k \alpha)|> |{\sin}(\pi k \re(\alpha))|> c k b^{-k},\qquad
k\ge1,
\]
which implies that the second series in (\ref{eq_lnphi_series}) can be
bounded from above by
\[
c^{-1} e^{\pi\im(\alpha)} \sum_{k\ge1} k^{-1} \bigl(b|z|^{\re
(\alpha)} \bigr)^k.
\]
It is clear that the above series converges in the domain $|z|^{\re
(\alpha)}<b$, $|{\arg}(z)|<\pi$. Now we take the limit $\im(\alpha)
\to
0^+$ and use the Dominated Convergence theorem to conclude that (\ref
{eq_lnphi_series}) is true when $z$ belongs to the above mentioned
domain and $\alpha$ is real,
$\alpha\notin{\mathcal L}\cup\q$.
Since $b>1$ is arbitrary we conclude that the series in (\ref
{eq_lnphi_series}) converges in the domain $|z|<1$, $|{\arg}(z)|<\pi$.
\end{pf}

Theorem \ref{thm_lnphi_series} was established independently in a
recent paper by Graczyk and Jakubowski \cite{Graczyk2009}, where the
authors applied series expansion to
the integrand in the formula for $\frac{\dd}{\dd z} \ln(\phi(z))$, the
latter was obtained from the Darling's integral.

We see that convergence of series (\ref{eq_lnphi_series}) is intimately
linked with the degree to which we can approximate a number by rational
numbers. The good news is that the series converges for $|z|<1$ for
almost all values of $\alpha>0$, as the set ${\mathcal L}$ has
Lebesgue measure zero. However, we can easily exhibit a dense
set of irrational numbers $\alpha$ for
which the argument that we used to prove this theorem breaks down. For
example, let us take an integer $a>1$ and define
\[
\alpha=\sum_{m\ge0}\frac{1}{q_m}\qquad \mbox{where }
q_{n+1}=a^{q_n}, q_0=1.
\]
It is clear that if we take the sum of the first $n$ terms we will have
a rational approximation in the form $p_n/q_n$, and the error of this
approximation will be less then $C a^{-q_{n}}$ for some $C>0$. Then,
assuming that $\rho$ is not a Liouville number and considering a
subseries of (\ref{eq_lnphi_series}) corresponding
to $k=q_m$ we obtain
\[
\biggl|\sum_{m\ge1} \frac{\sin(\pi q_m \alpha\rho)}{q_m\sin
(\pi q_m \alpha)}
(-1)^m z^{\alpha q_m} \biggr| > C_2 \sum_{m\ge1} \frac
{(az^{\alpha})^{q_m}}{p_m^{\mu} q_m},
\]
and we see that this series cannot converge unless $|az^{\alpha}|<1$,
therefore
the domain of convergence is strictly smaller then $|z|<1$. It is clear
that we will have the same situation if we multiply $\alpha$
by any rational number, thus we have a dense set of real parameters
$\alpha$ for which the domain of convergence of the series in (\ref
{eq_lnphi_series}) can be arbitrarily small. Of course the Lebesgue
measure of the set of such ``unlucky'' values of $\alpha$ is zero, so
Theorem \ref{thm_lnphi_series} can still be used (with some caution)
for numerical computations.

When $\im(\alpha)>0$, one can obtain another series representation for
$\phi(z)$ by applying the q-binomial theorem
(see Theorem 10.2.1 in \cite{Andrews})
\[
\frac{(az;q)_{\infty}}{(z;q)_{\infty}}=\sum_{n\ge0} \frac
{(a;q)_n}{(q;q)_n}z^n, \qquad |z|<1, |q|<1
\]
to (\ref{eq_phi_q_product}). However, it seems to be very hard to say
anything about the convergence of such series for $\alpha\in\rr
^{\plus}$
and we did not pursue this further.

In Corollary \ref{corollary_phi_analytic}, we have seen that $\phi
(e^w)$ is quasi-periodic with two different periods.
In the next theorem we collect other results on functional equations
satisfied by $\phi(z)$.
\begin{theorem}\label{thm_phi_func_eqns} Assume that $\rho\in(0,1)$.

\begin{enumerate}[(iii)]
\item[(i)] For $\alpha>0$ and $|{\arg}(z)|<\pi$
%
%
\begin{equation}\label{functional_eq_z_inv}
\phi\biggl({\frac{1}z};\alpha,\rho\biggr)= z^{\alpha\rho} \phi(z;\alpha
,\rho).
\end{equation}

\item[(ii)] For $\alpha\in(0,\rho^{-1})$ and $|{\arg}(z)|<\pi\min\{1,
\frac{1}{\alpha} \}$
%
%
\begin{equation}\label{functional_eq_z_alpha}
\phi(z;\alpha,\rho)=\phi\biggl(z^{\alpha}; \frac{1}{\alpha},\alpha
\rho\biggr).
\end{equation}

\item[(iii)] For $n \ge2$, $\alpha> n-1$ and $|{\arg}(z)| < \pi
(1-\frac{n-1}{\alpha} )$
%
%
\begin{equation}\label{multiplication_formula}
\phi(z; \alpha, \rho)=\prod_{k=0}^{n-1}
\phi\biggl(ze^{ (n-2k-1){\pi\ii}/{\alpha}}; \frac{\alpha}n,\rho\biggr).
\end{equation}
\item[(iv)] For $\alpha>0$ and $\im(z)\ne0$
%
%
\begin{equation}\label{functional_eq_WH_fact}
\phi(z;\alpha,1-\rho)\phi(-z;\alpha,\rho)=(1+e^{-\delta\pi\ii
\alpha\rho}z^{\alpha} )^{-1},
\end{equation}
where $\delta=\operatorname{sign}(\im(z))$.
\end{enumerate}
\end{theorem}
\begin{pf}
Statement (i) follows from Theorem \ref{thm_phi_Gamma2} and statement
(ii) from the series representation (\ref{eq_lnphi_series}).
To prove (iii) we use
Proposition \ref{prop_main} and Theorem \ref{thm_props_F_z_tau}, (ii),
to find that for $z\in\rr$
\[
\phi(z; \alpha, \rho)=\prod_{k=0}^{n-1}
\phi\biggl(z^{{1/n}} e^{ (n-2k-1){\pi\ii}/{n}}; \alpha n,
\frac{\rho}{n} \biggr).
\]
To obtain (\ref{multiplication_formula}) we apply transformation (\ref
{functional_eq_z_alpha}) in both parts of the above identity and
rescale parameters.
To prove statement (iv) we note that for $(\alpha,\rho)\in\aaa$
functions $\phi(-\ii z; \alpha, 1-\rho)$ and $\phi(\ii z;\alpha
,\rho)$ are
the Wiener--Hopf factors $\phi^{\plus}_1(z)$ and $\phi^{\minus}_1(z)$
thus functional equation
(\ref{functional_eq_WH_fact}) is just another way of writing the
Wiener--Hopf factorization $\phi^{\plus}_1(z)\phi^{\minus
}_1(z)=(1+\Psi
(z))^{-1}$.
The result for general $\alpha>0$ and $\rho\in(0,1)$ follows by
analytic continuation.
\end{pf}

Identity (\ref{functional_eq_z_inv}) first appeared in \cite{Fourati}
and later in \cite{Graczyk2009}, while (\ref{functional_eq_z_alpha}) is
equivalent to Theorem 3 in \cite{Doney1987} (it also appears in \cite
{Fourati}). Note that both statements (i) and (ii) of Theorem \ref
{thm_phi_func_eqns} can be established directly with the help of
Proposition \ref{prop_main} and results presented in Theorem \ref
{thm_props_F_z_tau}, (i) and (iv). Another
possible approach to establish these results is to use the same
technique as in the proof of Theorem 3 in \cite{Doney1987}: first we
verify the transformation identity for the processes
belonging to one of the $\cc_{k,l}$ classes and then prove the general
case by using the fact that classes $\cc_{k,l}$ are dense in the family
of stable processes. Note that Theorem \ref{thm_lnphi_series} and
functional equation (\ref{functional_eq_z_inv}) provide a convenient
method to compute the
Wiener--Hopf factors, as long as $|z|\ne1$ and $\alpha$ is not too
close to being a rational number.

\section{Mellin transform of the supremum functional}\label
{section_mellin_transform}

In the case of a general L\'evy process the Wiener--Hopf factor $\phiqp
(z)$ gives us\vspace*{1pt} valuable information about the distribution
of $S_{\ee(q)}$, where $\ee(q) \sim \operatorname{Exp}(q)$, but the only way to
translate this into information about the distribution of $S_t$ is
to perform a numerical Laplace or Fourier transform inversion in the
$q$-variable. The scaling property of stable processes allows us to
replace the
integral transform in the $q$-variable with an integral transform in
the $z$-variable, and a surprising fact is that the latter can be
evaluated in closed form.
The main tool in this section will be the Mellin transform of the
supremum $S_1$, defined as
\[
\mm(w)=\mm(w;\alpha,\rho)=\e[ (S_1)^{w-1} ].
\]
The function $\mm(s)$ is well defined if $\re(s)$ is sufficiently close
to $1$. Note that our $\mm(s)$ corresponds to $\mm^{\plus}(s)$ in
\cite{Doney1987}. Next, we define the Mellin transfom of $\phi(z)$ as
\[
\Phi(s;\alpha,\rho)=\int_0^{\infty} z^{s-1} \phi(z;\alpha,\rho)
\,\dd z=\int_{\rr} e^{ws} \phi(e^w;\alpha,\rho) \,\dd w.
\]
The link between $\mm(s)$ and $\Phi(s)$ is the following identity,
which goes back to \cite{Darling1956} (see also \cite{Doney1987}):
%
%
\begin{equation}\label{identity_Phi_M}
\Phi(s;\alpha,\rho)=\Gamma(s) \Gamma\biggl( 1-\frac{s}{\alpha}\biggr)\mm
(1-s;\alpha,\rho).
\end{equation}
This identity can be easily established using the scaling property of
stable processes, which implies that $S_t\stackrel{d}{=} t^{
{1}/{\alpha}}S_1$
\begin{eqnarray*}
\Phi(s)&=&\int_0^{\infty} z^{s-1} \e[e^{-zS_{\ee(1)}}
] \,\dd z=
\Gamma(s) \e\bigl[ \bigl(S_{\ee(1)}\bigr)^{-s}\bigr]=\Gamma(s) \e
[ \ee(1)^{-{s}/{\alpha}} (S_1)^{-s} ]\\
&=&
\Gamma(s) \e[ \ee(1)^{-{s}/{\alpha}}] \e[
(S_1)^{-s} ]=
\Gamma(s) \Gamma\biggl( 1-\frac{s}{\alpha}\biggr)\mm(1-s;\alpha,\rho),
\end{eqnarray*}
where we have also used the fact that $\ee(1) \sim \operatorname{Exp}(1)$ is
independent of $S_t$.

With the help of formula (\ref{identity_Phi_M}) we can translate
functional equations for the Wiener--Hopf factor $\phi(z)$ into the
following remarkable identities for the Mellin transform $\mm(s)$.
\begin{theorem}\label{thm_Ms_func_eqns} The function $\mm(s)$ can be
analytically continued to a meromorphic function in $\ccc$.
$\mm(s)$ is quasiperiodic with periods $1$ and $\alpha$
%
%
\begin{eqnarray}\qquad
\label{eqn_M_period_1}
\mm(s+1)&=& \frac{\alpha}{\pi}\sin\biggl( \pi\biggl(\rho-\frac
{1-s}{\alpha}\biggr) \biggr) \Gamma\biggl(1-\frac{s}{\alpha}\biggr)\Gamma
\biggl(1-\frac{1-s}{\alpha}\biggr) \mm(s),
\\
\label{eqn_M_period_alpha}
\mm(s+\alpha)&=&\frac{\alpha}{\pi} \sin\bigl(\pi(\alpha\rho-1 +s)\bigr)
\Gamma
(1-s)\Gamma(\alpha-1+s) \mm(s),
\end{eqnarray}
and it satisfies
%
%
\begin{eqnarray}
\label{eq_M_trans1}
\mm(s;\alpha,\rho)&=&\frac{\Gamma(\alpha\rho-1+s)}{\Gamma(1-s)}
\frac{\Gamma(1-\rho+({1-s})/{\alpha} )}{\Gamma(1-
({1-s})/{\alpha})}\mm(2-\alpha\rho-s;\alpha,\rho), \hspace*{-34pt}\\
\label{eq_M_trans2}
\mm(s;\alpha,\rho)&=& \frac{\Gamma(s)}{\Gamma(2-s)}
\frac{\Gamma(1+({1-s})/{\alpha})}{\Gamma(1-
({1-s})/{\alpha})}
\mm\biggl(1-\frac{1-s}{\alpha};\frac{1}{\alpha
},\alpha\rho\biggr).
\end{eqnarray}
\end{theorem}
\begin{pf}
First of all, we use (\ref{functional_eq_z_inv}) to find that $\phi(z)
\sim z^{-\alpha\rho}$ as $z\to+\infty$. Thus the function
$\Phi(s;\alpha,\rho)$ is well defined for $0<\re(s)<\alpha\rho$.
Let us assume that $\alpha<\frac13$ and define
%
%
\begin{equation}\label{thm_MS_def_G}
G(s)=\int_{\rr} \frac{1-e^{-2\pi\ii\alpha\rho}}{1+e^{\alpha
w+\pi\ii\alpha(1+\rho)}} e^{ws} \phi(e^w) \,\dd w.
\end{equation}
Using the first functional equation in (\ref{phi_analytic_ctnn}) we
find the following identity:
%
%
\begin{eqnarray}\label{thm_Ms_proof_def_G}
G(s)&=&
\int_{\rr} e^{ws} \phi(e^{w+2\pi\ii}) \,\dd w-e^{-2\pi\ii\alpha
\rho} \int_{\rr} e^{ws} \phi(e^{w}) \,\dd w \nonumber\\[-8pt]\\[-8pt]
&=&
\int_{\rr} e^{ws} \phi(e^{w+2\pi\ii}) \,\dd w-e^{-2\pi\ii\alpha
\rho} \Phi(s).\nonumber
\end{eqnarray}
Corollary \ref{corollary_phi_analytic} and the condition $\alpha
<\frac
13$ guarantee that $\phi(e^w)$ is analytic for
$0<\im(w)<3\pi$; thus we can shift the contour of integration $\rr
\mapsto\rr+2\pi\ii$ in the integral in the right-hand side of the
above identity
and obtain
\[
\int_{\rr} e^{ws} \phi(e^{w+2\pi\ii}) \,\dd w=e^{-2\pi\ii s} \int
_{\rr} e^{ws} \phi(e^{w}) \,\dd w=e^{-2\pi\ii s} \Phi(s).
\]
Combining (\ref{thm_Ms_proof_def_G}) with the above equation we
conclude that
%
%
\begin{equation}\label{eqn_Gs_Phi}
G(s)=(e^{-2\pi\ii s}-e^{-2\pi\ii\alpha\rho}) \Phi(s).
\end{equation}
The integral which defines function $G(s)$ converges for $0<\re
(s)<\alpha(1+\rho)$,
and therefore (\ref{eqn_Gs_Phi}) allows us to continue the function
$\Phi(s)$ analytically into a larger domain $0<\re(s)<\alpha(1+\rho)$.
Next, one can use definitions of $\Phi(s)$ and $G(s)$ to check that
\[
G(s)+e^{\pi\ii\alpha(1+\rho)} G(s+\alpha)=(1-e^{-2\pi\ii\alpha
\rho}) \Phi(s),
\]
and combining the above two equations we obtain
\[
\Phi(s+\alpha)=-\frac{\sin(\pi s)}{\sin(\pi(s+\alpha(1-\rho)))}
\Phi(s).
\]
Formula (\ref{eqn_M_period_alpha}) follows easily by applying (\ref
{identity_Phi_M}) to the above identity.
Equation (\ref{eqn_Gs_Phi}) extends $\Phi(s)$ analytically into the
domain $0<\re(s)<\alpha(1+\rho)$, which means that
$\mm(s)$ is well defined for $1-\alpha(1+\rho)<\re(w)<1$ and applying
(\ref{eqn_M_period_alpha})
repeatedly we can extend $\mm(s)$ to a meromorphic function in the
entire complex plane. The restriction $\alpha<\frac13$
can be removed by analytic continuation to other values of $\alpha$ and
the second identity (\ref{eqn_M_period_1}) can be proved similarly by
starting with the second functional equation in (\ref{phi_analytic_ctnn}).

To prove (\ref{eq_M_trans1}) and (\ref{eq_M_trans2}), we apply the
Mellin transform to both sides of equations (\ref{functional_eq_z_inv})
and (\ref{functional_eq_z_alpha}) to obtain
\[
\Phi(s;\alpha,\rho)=\Phi(\alpha\rho-s;\alpha,\rho),\qquad \Phi
(s;\alpha
,\rho)=\frac{1}{\alpha}\Phi\biggl(\frac{s}{\alpha};\frac{1}{\alpha
},\alpha\rho\biggr).
\]
Functional equations (\ref{eq_M_trans1}) and (\ref{eq_M_trans2}) follow
from the above identities and (\ref{identity_Phi_M}), the details are
left to the reader.
\end{pf}

A surprising fact is that quasi-periodicity of $\mm(s)$ allows us to
find an expicit formula for this function, and having
already developed
all the necessary tools we can enjoy a rather simple and
straightforward proof.
\begin{theorem}\label{thm_M_explicit}
For $s\in\ccc$
%
%
\begin{eqnarray}\label{formula_Ms_general}
\mm(s)&=&\alpha^{s-1} \frac{G(\alpha\rho;\alpha)}{G(\alpha
(1-\rho
)+1;\alpha)}
\nonumber\\[-8pt]\\[-8pt]
&&{}\times \frac{G(\alpha(1-\rho)+2-s;\alpha)}
{G(\alpha\rho-1+s;\alpha)}
\times\frac{G(\alpha-1+s;\alpha)}{G(\alpha+1-s;\alpha)}.\nonumber
\end{eqnarray}
If $X_t\in\cc_{k,l}$ and $l>0$, then
%
%
\begin{eqnarray}\label{formula_Ms_Cklplus}
\mm(s)&=&
\frac{\Gamma(s)}{\Gamma(1-({1-s})/{\alpha})}
\prod_{j=1}^{l-1} \frac{\sin({\pi}/{\alpha}(s-1+j)
)}{\sin({\pi j}/{\alpha})}\nonumber\\[-8pt]\\[-8pt]
&&{}\times
\prod_{j=1}^{k} \frac{\sin(\pi\alpha j)}{\sin(\pi(1-s+\alpha
j))}.\nonumber
\end{eqnarray}
If $X_t\in\cc_{k,l}$ and $l<0$, then
%
%
\begin{eqnarray}\label{formula_Ms_Cklminus}
\mm(s)&=&\frac{\Gamma(1+({1-s})/{\alpha})}{\Gamma(2-s)}
\prod_{j=1}^{|k|-1} \frac{ \sin(\pi(s-1+\alpha j))}{\sin(\pi
\alpha j) }\nonumber\\[-8pt]\\[-8pt]
&&{}\times\prod_{j=1}^{|l|} \frac{\sin({\pi j}/{\alpha}
)}{\sin({\pi}/{\alpha}(1-s+j))}.\nonumber
\end{eqnarray}
\end{theorem}
\begin{pf}
Formula (\ref{formula_Ms_general}) can be established using the
classical approach from the theory of elliptic functions. Let us
denote the function on the right-hand side of this
formula as $H(s)$. Using quasi-periodicity properties of $G(z;\tau)$
[see formulas (\ref{funct_rel_G})] we find that $H(s)$ satisfies
(\ref{eqn_M_period_1}) and (\ref{eqn_M_period_alpha}) (the detailed
computations can be found in the Appendix of the online version of this
paper, see \href{http://arxiv.org/abs/1001.0991}{arXiv:1001.0991}).
This implies that $\mm(s)/H(s)$ is a meromorphic function which is
periodic with
periods $1$ and $\alpha$. An integer linear combination of
periods is also a period; thus, assuming that $\alpha\in\rr^{\plus}$
is irrational we can find arbitrarily small periods of the form
$m\alpha+n$, and therefore the function $\mm(s)/H(s)$ must be constant.
The value of this constant is equal to one since $\mm(1)=H(1)=1$.
Formulas (\ref{formula_Ms_Cklplus}) and (\ref{formula_Ms_Cklminus}) can
be verified using exactly the same approach, again the detailed proof
can be found in the Appendix of the online version of this paper, see
\href{http://arxiv.org/abs/1001.0991}{arXiv:1001.0991}. Otherwise,
formulas (\ref{formula_Ms_Cklplus}) and (\ref{formula_Ms_Cklminus})
can be established directly, using the fact that
$\alpha\rho=l-k\alpha$ and transforming (\ref{formula_Ms_general})
with the help of the functional equations for the double gamma function
(\ref{funct_rel_G}). Finally, we extend the result to the case when
$\alpha$ is rational by analytic continuation (note that all the formulas
in Theorem \ref{thm_M_explicit} are well defined when $\alpha\in\q$).
\end{pf}

The Mellin transform of the supremum was evaluated for the three
special cases: $\cc_{0,1}$ (spectrally negative process) and $\cc
_{1,1}$, $\cc_{1,2}$ (processes with two-sided jumps) in Corollary 2
\cite{Doney1987}, but they do not agree with our results presented in
Theorem \ref{thm_M_explicit}.
It seems that formulas (3.4), (3.5) and (3.6) in \cite{Doney1987} are
not correct, as they violate identity $\mm(1)=1$ (which follows from
the definition of the Mellin transform).

It is interesting to note that formula (\ref{formula_Ms_general}) is
valid even in the limiting case $\alpha=2$ and $\rho=1/2$, that is,
$X_t$ is Brownian motion $\sqrt{2} W_t$. Using the properties of the
double gamma function given in (\ref{funct_rel_G}) one can check that
(\ref{formula_Ms_general})
reduces to
\[
\mm(s)=2^{s-1}\frac{\Gamma({s}/{2})}{\Gamma(1/2
)},
\]
and the above formula can also be verified directly using the
distribution of the supremum of Brownian motion (see \cite{Borodin}).

The results in Theorems \ref{thm_Ms_func_eqns} and \ref{thm_M_explicit}
are surprisingly explicit, and later in Section \ref
{section_pdf_supremum} we
will use them to derive series representations and asymptotic
expansions for the probability density function of $S_1$.
Here we also present two immediate corollaries of Theorems \ref
{thm_Ms_func_eqns} and \ref{thm_M_explicit}.
\begin{corollary}\label{corollary_id1} Assume that\vspace*{1pt} $\alpha\in(\frac
12,1) \cup(1,2)$. Let $X_t$ and $\tilde X_t$ be stable processes with
parameters $(\alpha,\rho)$ and $(\frac{1}{\alpha},\alpha\rho)$ having
corresponding supremum processes $S_t$ and $\tilde S_t$. Then we have
the following identity in distribution:
%
%
\begin{equation}\label{eqn_distr_id1}
\varepsilon_1\biggl[\frac{S_1}{\varepsilon_2}\biggr]^{\alpha} \stackrel
{d}{=} \varepsilon_3^{\alpha}\biggl[\frac{\tilde S_1}{\varepsilon_4}
\biggr],
\end{equation}
where $\varepsilon_i \sim \operatorname{Exp}(1)$ and all random variables are assumed
to be independent.
\end{corollary}
\begin{pf}
To prove this statement we change the variable in (\ref{eq_M_trans2})
$s \mapsto1-\alpha(1-s)$, rewrite the resulting identity as:
\[
\Gamma(1+\alpha-\alpha s)\Gamma(s)\mm(1-\alpha+\alpha s;\alpha
,\rho)=
\Gamma(1-\alpha+\alpha s) \Gamma(2-s)
\mm\biggl(s;\frac{1}{\alpha},\alpha\rho\biggr)
\]
and use the following facts: (i) $\Gamma(s)$ is the Mellin transform of
the $\operatorname{Exp}(1)$ random variable; (ii)
if $M(s)=\e[X^{s-1}]$ is the Mellin transform of the
random variable $X$, then $M(1-a+as)$ is the Mellin transform of
$Y=X^a$; (iii) the Mellin transform of the
product of two independent random variables is equal to the product of
their Mellin transforms.
\end{pf}

Corollary \ref{corollary_id1} should be compared with Theorems 12 and
13 in \cite{Zolotarev1957}, which describe
identities in distribution satisfied by products of powers of stable
random variables, and with Theorems 1.1, 1.2 in \cite{Fourati}.
Also, note that using the scaling property of stable processes
we can rewrite identity (\ref{eqn_distr_id1}) in a more symmetric form
\[
\biggl[\frac{S_{\ee(1)}}{\varepsilon_1}\biggr]^{\alpha} \stackrel{d}{=}
\frac{\tilde S_{\ee(1)}}{\varepsilon_2} .
\]
\begin{corollary}\label{corollary_id2} If $X_t \in\cc_{k,l}$ and
$l>0$, then $S_1$ satisfies the following identity in distribution:
%
%
\begin{equation}\label{eqn_distr_id2}
S_1 \times\Biggl[ \varepsilon_1
\prod_{j=1}^{l-1} \frac{\gamma_{\{\alpha^{-1} j \}}}{\gamma_{1-\{
\alpha^{-1} j \}}} \Biggr]^{1/\alpha}
\stackrel{d}{=}
\varepsilon_2
\prod_{j=1}^{k} \frac{\gamma_{1-\{\alpha j \}}}{\gamma_{\{\alpha
j \}}} ,
\end{equation}
where $\{x\} \in[0,1) $ denotes the fractional part of $x$,
$\varepsilon_i \sim \operatorname{Exp}(1)$, $\gamma_k$ denotes a gamma random
variable with
$\e[\gamma_k]=\operatorname{Var}[\gamma_k]=k$ and all random variables
are assumed to be independent.
\end{corollary}
\begin{pf}
The proof is very similar to the proof of Corollary \ref
{corollary_id1}: we rewrite (\ref{formula_Ms_Cklplus}) as
\[
\mm(s) \Gamma\biggl(1-\frac{1-s}{\alpha}\biggr)
\prod_{j=1}^{l-1} \frac{\sin({\pi j}/{\alpha}
)}{\sin({\pi}/{\alpha}(s-1+j))}=\Gamma(s)
\prod_{j=1}^{k} \frac{\sin(\pi\alpha j)}{\sin(\pi(1-s+\alpha j))},
\]
and use the fact that for $0<a<1$ the Mellin transform of $Y=\gamma
_a/\gamma_{1-a}$ (where
$\gamma_a$ and $\gamma_{1-a}$ are independent) can be computed as
\begin{eqnarray*}
\e[ Y^{s-1} ]&=&\e[ (\gamma_a)^{s-1}
]\e[ (\gamma_{1-a})^{1-s} ] =
\frac{\Gamma(a+s-1)}{\Gamma(a)}\frac{\Gamma(2-a-s)}{\Gamma(1-a)}\\[-7pt]
&=&
\frac{\sin(\pi a)}{\sin(\pi(a+s-1))}.
\end{eqnarray*}
\upqed
\end{pf}

Corollary \ref{corollary_id2} should be considered as a generalization
of the corresponding result for Brownian motion, that is, when
$X_t=\sqrt{2} W_t$.
In this case we have
\[
S_1 \times\varepsilon_1^{1/2}\stackrel{d}{=}\varepsilon_2.
\]
This fact can be established by applying the scaling property to
the identity $S_{\ee(1)} \sim \operatorname{Exp}(1)$, which follows from the
Wiener--Hopf factorization for Brownian motion. Identity (\ref
{eqn_distr_id2}) also can be rewritten in the form
$S_{T_1}=T_2$ with an obvious choice of random variables $T_i$. It
might be possible to establish a similar identity in the general case,
using formula (\ref{formula_Ms_general}) instead of (\ref
{formula_Ms_Cklplus}). However, it is not clear how to connect the
double gamma function with the
Mellin transform of some random variable and we will leave this for
future work.

\section{Probability density function of the supremum
functional}\label
{section_pdf_supremum}

In the last two years there have appeared several interesting and
important results related to the density of the supremum $S_1=\sup\{
X_s\dvtx  0\le s\le1\}$.
In the spectrally positive case Doney \cite{Doney2008} has obtained the
first asymptotic term for the density
\[
p(x)=\frac{\dd}{\dd x} \p(S_1 \le x)
\]
as $x\to+\infty$; Bernyk, Dalang and Peskir \cite{Bernyk2008} have
derived an explicit convergent series representation for $p(x)$ and
Patie \cite{Patie2009}
has obtained a complete asymptotic expansion of $p(x)$ as $x\to+\infty$.
In the case of a general stable process Doney and Savov \cite{Doney2010} obtain
the first term of asymptotic expansion of $p(x)$ as $x \to0^+$ or
$x\to
+\infty$, and they also mention that it is possible to
obtain higher order asymptotic terms as $x\to+\infty$ (see Section 5
in their paper). See also \cite{Vatutin2009} for asymptotic results
on distributions of functionals of a certain random walk related to
stable processes and \cite{Peskir2008} for an explicit infinite series
representation for the density of the first hitting time of a point in
the spectrally positive case.

The following two theorems summarize our main results in this section.
\begin{theorem}\label{thm1} Assume that $\alpha\notin\q$. Define
sequences $\{a_{m,n}\}_{m\ge0,n\ge0}$ and $\{b_{m,n}\}_{m\ge0,n\ge
1}$ as
%
%
\begin{eqnarray}\label{def_a_mn}
a_{m,n}&=&\frac{(-1)^{m+n} }{\Gamma(1-\rho-n-{m}/{\alpha}
)\Gamma(\alpha\rho+m+\alpha n)} \nonumber\\[-8pt]\\[-8pt]
&&{}\times
\prod_{j=1}^{m} \frac{\sin({\pi}/{\alpha} ( \alpha
\rho+ j-1 ))} {\sin({\pi j}/{\alpha} )}
\prod_{j=1}^{n} \frac{\sin(\pi\alpha(\rho+j-1))}{\sin(\pi
\alpha j)},\nonumber
\\
%
%
\label{def_b_mn}
b_{m,n}&=&\frac{\Gamma(1-\rho-n-{m}/{\alpha})\Gamma(\alpha
\rho+m+\alpha n) }{\Gamma(1+n+{m}/{\alpha})\Gamma
(-m-\alpha n)}
a_{m,n}.
\end{eqnarray}
Then we have the following asymptotic expansions:
%
%
\begin{eqnarray}
\label{eqn_p_0}
p(x) &\sim& x^{\alpha\rho-1} \sum_{n\ge0} \sum_{m\ge0}
a_{m,n} x^{ m+\alpha n},\qquad x\to0^+, \\
\label{eqn_p_infty}
p(x) &\sim& x^{-1-\alpha} \sum_{n\ge0} \sum_{m\ge
0}b_{m,n+1} x^{-m-\alpha n},\qquad x\to+\infty.
\end{eqnarray}
\end{theorem}
\begin{remark}
Asymptotic expansion (\ref{eqn_p_0}) (and other similar expressions)
should be understood in the following sense: for any $c>0$ we have
\[
p(x) = x^{\alpha\rho-1} \mathop{\sum_{m,n\ge0}}_{m+\alpha n < c
}a_{m,n} x^{ m+\alpha n}
+O(x^{\alpha\rho-1+c}),\qquad x\to0^+.
\]
\end{remark}
\begin{theorem}\label{thm2} Assume that $X_t \in\cc_{k,l}$.
If $l>0$, then for $n\in\{0,1,\ldots,k\}$ and $m \in\zz$, we define
%
%
\begin{eqnarray}\label{def_c_plus}
c^+_{m,n}&=& \frac{(-1)^{m(k+1)+nl+1} }{\Gamma(1+n+{m}/{\alpha
})\Gamma(-m-\alpha n)} \nonumber\\[-8pt]\\[-8pt]
&&{}\times
\prod_{j=1}^{l-1} \frac{\sin({\pi}/{\alpha}(j+m)
)} {\sin({\pi j}/{\alpha} )}
\prod_{j=1}^{k-n} \frac{\sin(\pi\alpha(j+n))}{\sin(\pi\alpha j)},\nonumber
\end{eqnarray}
while if $l<0$, then for $m \in\{0,1,\ldots,|l|\}$ and $n \in\zz$, we define
%
%
\begin{eqnarray}\label{def_c_minus}
c^-_{m,n}&=& \frac{(-1)^{mk+n(l+1)+1} }{\Gamma(1+n+{m}/{\alpha
})\Gamma(-m-\alpha n)}\nonumber\\[-8pt]\\[-8pt]
&&{}\times
\prod_{j=1}^{|k|-1} \frac{\sin(\pi\alpha(j+n))} {\sin(\pi
\alpha j)}
\prod_{j=1}^{|l|-m} \frac{\sin({\pi}/{ \alpha}
(j+m))}{\sin({\pi j}/{\alpha} )}. \nonumber
\end{eqnarray}
Then if $\alpha\in(0,1)$ and $l>0$ $\{l<0\}$ we have a convergent
series representation
%
%
\begin{eqnarray}
\label{alpha_01_eq_1}
p(x) &=& -\sum_{n=1}^k \sum_{m \ge0} c^+_{m,n}
x^{-m-\alpha n-1}, \nonumber\\[-8pt]\\[-8pt]
\Biggl\{ p(x) &=& -\sum_{m=0}^{|l|} \sum_{n \ge1}
c^-_{m,n} x^{-m-\alpha n-1} \Biggr\},\qquad x \in\rr^+, \nonumber
\end{eqnarray}
and an asymptotic expansion
%
%
\begin{eqnarray}
\label{alpha_01_eq_2}
p(x) &\sim& \sum_{n=0}^k \sum_{m \le-l} c^+_{m,n}
x^{-m-\alpha n-1}, \nonumber\\[-8pt]\\[-8pt]
\Biggl\{ p(x) &\sim& \sum_{m=1}^{|l|} \sum_{n \le k}
c^-_{m,n} x^{-m-\alpha n-1} \Biggr\},\qquad x \to0^+ . \nonumber
\end{eqnarray}
Similarly, if $\alpha\in(1,2)$ and $l>0$ $\{l<0\}$,
we have a convergent series representation
%
%
\begin{eqnarray}\label{alpha_12_eq_1}
p(x) &=& \sum_{n=0}^k \sum_{m \le-l} c^+_{m,n}
x^{-m-\alpha n-1}, \nonumber\\[-8pt]\\[-8pt]
\Biggl\{ p(x) &=& \sum_{m=1}^{|l|} \sum_{n \le k}
c^-_{m,n} x^{-m-\alpha n-1} \Biggr\},\qquad x \in\rr^+, \nonumber
\end{eqnarray}
and an asymptotic expansion
%
%
\begin{eqnarray}\label{alpha_12_eq_2}
p(x) &\sim& -\sum_{n=1}^k \sum_{m \ge0} c^+_{m,n}
x^{-m-\alpha n-1}, \nonumber\\[-8pt]\\[-8pt]
\Biggl\{ p(x) &\sim& -\sum_{m=0}^{|l|} \sum_{n \ge1}
c^-_{m,n} x^{-m-\alpha n-1} \Biggr\},\qquad x \to+\infty. \nonumber
\end{eqnarray}
Infinite series in (\ref{alpha_01_eq_1}) and (\ref{alpha_12_eq_1})
converge uniformly on compact subsets of $\rr^+$.
\end{theorem}

First we need to establish some technical results which describe the
analytic structure of $\mm(s)$.
\begin{lemma}\label{lemma1} The function $\mm(s)$ can be analytically
continued to a meromorphic function in $\ccc$. If $\alpha\notin\q$
and $X_t \notin\cc_{k,l}$ for all $k$ and $l$, then $\mm(s)$ has
simple poles at
\begin{eqnarray*}
\{s^+_{m,n}\}_{m\ge1, n\ge1}&=&\{m+\alpha n\}_{m\ge1, n\ge1}, \\
\{s^-_{m,n}\}_{m\ge0, n\ge0}&=&\{1-\alpha\rho-m- \alpha n\}_{m\ge
0, n\ge0},
\end{eqnarray*}
with residues
%
%
\begin{equation}\label{res_amn_bmn}
\operatorname{Res}\bigl(\mm(s) \dvtx  s^+_{m,n} \bigr)=-b_{m-1,n},\qquad
\operatorname{Res}\bigl(\mm(s) \dvtx  s^-_{m,n} \bigr)=a_{m,n}.
\end{equation}
In the case when $X_t\in\cc_{k,l}$ and $l>0$ $\{l<0\}$, the function
$\mm(s)$ has simple poles at $s_{m,n}=m+\alpha n$, where
\begin{eqnarray*}
m &\le& 1-l,\qquad n\in\{0,1,2,\ldots,k\} \quad\mbox{or}\quad
m \ge1,\qquad n\in\{1,2,\ldots,k\}, \\
\bigl\{ m&\in&\{1,2,3,\ldots,|l|+1\}, n \ge1 \mbox{ or }
m\in\{2,3,\ldots,|l|+1\}, n \le k \bigr\},
\end{eqnarray*}
with residues
%
%
\begin{equation}\label{res_cmn}
\operatorname{Res}\bigl(\mm(s) \dvtx  s_{m,n} \bigr)=c^+_{m-1,n},\qquad
\bigl\{ \operatorname{Res}\bigl(\mm(s) \dvtx  s_{m,n} \bigr)=c^-_{m-1,n}
\bigr\}.
\end{equation}
\end{lemma}
\begin{pf} Iterating identities (\ref{eqn_M_period_1}) and (\ref
{eqn_M_period_alpha}), we find
\begin{eqnarray*}
&&
\mm(1+s+m+\alpha n)\\
&&\qquad=\mm(1+s) \alpha^{m+n} \\
&&\qquad\quad{}\times\prod_{j=0}^{m-1}
\frac{\Gamma(1-({1+s+j})/{\alpha})\Gamma(1+
({s+j})/{\alpha})}
{\Gamma(1-\rho-({s+j})/{\alpha})\Gamma(\rho+
({s+j})/{\alpha})} \\
&&\qquad\quad{}\times
\prod_{j=0}^{n-1} \frac{\Gamma(-s-m-\alpha j)\Gamma(\alpha
+s+m+\alpha j)}
{\Gamma(1-\alpha\rho-s-m-\alpha j) \Gamma(\alpha\rho+s +m +
\alpha
j)} .
\end{eqnarray*}
The right-hand side of the above equation has a simple pole at $s=0$,
which comes from $\Gamma(-s-m)$ (take $j=0$ in the last product).
The residue of $\Gamma(-s-m)$ at $s=0$ is equal to $(-1)^{m-1}/m!$, and
therefore
%
%
\begin{eqnarray}\label{eq_res_mms}
&& \operatorname{Res}\bigl(\mm(s) \dvtx  1+m+\alpha n\bigr)\nonumber\\
&&\qquad= \frac{(-1)^{m-1}\alpha
^{m+n}\Gamma(\alpha+m)}{m! \Gamma(1-\alpha\rho-m) \Gamma(\alpha
\rho
+m)} \nonumber\\[-8pt]\\[-8pt]
&&\qquad\quad{}\times \prod_{j=0}^{m-1} \frac{\Gamma(1-
({1+j})/{\alpha})\Gamma(1+{j}/{\alpha})}
{\Gamma(1-\rho-{j}/{\alpha})\Gamma(\rho+
{j}/{\alpha})}\nonumber\\
&&\qquad\quad{}\times
\prod_{j=1}^{n-1} \frac{\Gamma(-m-\alpha j)\Gamma(\alpha
+m+\alpha j)}
{\Gamma(1-\alpha\rho-m-\alpha j) \Gamma(\alpha\rho+m + \alpha j)}.\nonumber
\end{eqnarray}
The above expression is equal to $-b_{m,n}$, which can be verified with
the help of the reflection formula for the gamma function. The detailed
computations
can be found in the Appendix of the online version of this paper, see
\href{http://arxiv.org/abs/1001.0991}{arXiv:1001.0991}.
The value of the residue at $s^-_{m,n}$ can be easily established with
the help of the reflection formula (\ref{eq_M_trans1}). Finally, when
$X_t \in\cc_{k,l}$ equation
(\ref{res_cmn}) can be derived by computing the residues of $\mm(s)$
given by formula (\ref{formula_Ms_Cklplus}). Otherwise one could establish
it as a corollary of the previous result by checking that
$a_{m,n}=c^{\pm}_{-l-m,k-n}$ and $b_{m,n}=-c^{\pm}_{m,n}$ depending on
$\pm l>0$.
The latter approach with detailed calculations can be found in the
Appendix of the online version of this paper, see \href
{http://arxiv.org/abs/1001.0991}{arXiv:1001.0991}.
\end{pf}
\begin{lemma}\label{lemma2}
For $x\in\rr$ we have as $y \to\infty$, $y\in\rr$
%
%
\begin{eqnarray}\label{asympt_M_i_infty}
\ln\bigl(|\mm(x+\ii y)|\bigr) &=& -\frac{\pi|y|}{2\alpha} \bigl( \alpha(1-\rho
)+1-\alpha\rho\bigr)\nonumber\\[-8pt]\\[-8pt]
&&{} +o(y).\nonumber
\end{eqnarray}
\end{lemma}
\begin{pf}
Equation (4.5) in \cite{Bill1997} gives us the following asymptotic
expansion valid as $z\to\infty$ in the domain $|{\arg}(z)|<\pi$:
%
%
\begin{eqnarray}\label{G_asympt_formula}
\ln( G(z;\alpha))&=&\frac{1}{2\alpha}z^2 \ln(z) - \frac{3+2\ln
(\alpha
)}{4\alpha} z^2 - \frac{1+\alpha}{2\alpha} z \ln(z)
\nonumber\\
&&{}+ \frac12 \biggl( \frac{1+\ln(\alpha)}{\alpha}+ \ln(2\pi
\alpha) +1\biggr)z
\\
&&{}+ \biggl(\frac{\alpha}{12}+\frac14+\frac{1}{12\alpha}\biggr)
\ln(z) + c(\alpha)+O(1/z),\nonumber
\end{eqnarray}
where $c(\alpha)$ is some constant depending on $\alpha$. The
asymptotic formula (\ref{asympt_M_i_infty}) can be obtained from
(\ref{formula_Ms_general}) and (\ref{G_asympt_formula}) using the
following asymptotic expansion for the logarithm function:
\[
\ln(A+s)=\ln(s)+\frac{A}{s}-\frac{A^2}{2s^2}+O(s^{-3}),\qquad s\to
\infty,\qquad |{\arg}(z)|<\pi,
\]
after some straightforward but tedious computations, which can be
considerably simplified with the help of symbolic computation software.
\end{pf}
\begin{pf*}{Proof of Theorem \ref{thm1}}
Equation (\ref{asympt_M_i_infty}) and the fact that $\alpha\rho\le1$
for all $(\alpha,\rho)\in\aaa$ imply
that $\mm(s)$ decreases exponentially as $\im(s) \to\infty$; thus
$S_1$ has a smooth density function $p(x)$,
which can be obtained as the inverse Mellin transform,
%
%
\begin{equation}\label{p_as_inv_of_m}
p(x)=\frac{1}{2\pi\ii} \int_{1+\ii\rr} \mm(s) x^{-s} \,\dd s.
\end{equation}
By shifting the contour of integration
$1+\ii\rr\mapsto c+\ii\rr$ where $c<0$ and taking care of the
residues at points $s^-_{m,n}$ we obtain
%
%
\begin{equation}\label{p_shift_contour}\quad
p(x)=\sum\operatorname{Res}\bigl(\mm(s)\dvtx  s^-_{m,n}\bigr) \times x^{-s^-_{m,n}}+
\frac{1}{2\pi\ii} \int_{c+\ii\rr} \mm(s) x^{-s} \,\dd s,
\end{equation}
where the summation is over all $m\ge0, n \ge0$, such that
$s^-_{m,n}>c$. The integral in the right-hand side of (\ref
{p_shift_contour}) can be estimated as follows:
\begin{eqnarray*}
\biggl| \frac{1}{2\pi\ii} \int_{c+\ii\rr} \mm(s) x^{-s} \,\dd s
\biggr| &=& \frac{x^{-c}}{2\pi}
\biggl| \int_{\rr} \mm(c+\ii t) x^{-it} \,\dd t \biggr| \\
&<&
\frac{x^{-c}}{2\pi}
\int_{\rr} |\mm(c+\ii t)| \,\dd t = O(x^{-c}).
\end{eqnarray*}
This establishes the asymptotic expansion (\ref{eqn_p_0}). The proof of
(\ref{eqn_p_infty}) is identical, except that now we need to shift the contour
of integration in (\ref{p_as_inv_of_m}) in the opposite direction.
\end{pf*}
\begin{pf*}{Proof of Theorem \ref{thm2}}
The asymptotic expansions (\ref{alpha_01_eq_2}) and (\ref
{alpha_12_eq_2}) can be derived in the same way as in the proof of
Theorem \ref{thm1} (or obtained
as its corollaries); thus we only need to establish convergence of
(\ref
{alpha_01_eq_1}) and (\ref{alpha_12_eq_1}).
Let us assume that $\alpha\in(1,2)$ and $l>0$. Choose $c_0 \in(0,1)
$ such that $c_0 \ne m+\alpha n $ for $n\in\{0,1,2,\ldots,k\}$ and $m \in
\zz$.
We take $N$ to be a large positive number and shift the contour of
integration in (\ref{p_as_inv_of_m})
$1+\ii\rr\mapsto c_0-N+\ii\rr$ while taking into account residues at
$s_{m,n}$, hence
%
%
\begin{equation}\label{p_shift_contour2}\qquad
p(x)=\sum\operatorname{Res}\bigl(\mm(s)\dvtx  s_{m,n}\bigr) \times x^{-s_{m,n}}+
\frac
{1}{2\pi\ii} \int_{c_0-N+\ii\rr} \mm(s) x^{-s} \,\dd s,
\end{equation}
where the summation is over $m,n$ such that $c_0-N<\re(s_{m,n})<1$.

Using (\ref{formula_Ms_Cklplus}) and the reflection formula
for the gamma function we find that for some constant $C \in\rr$
%
%
\begin{equation}\label{mm_Ckl}
\mm(s)= C \frac{\Gamma(({1-s})/{\alpha})} {\Gamma(1-s)}
\frac{
\prod_{j=0}^{l-1} \sin({\pi}/{\alpha}(s-1+j))}{
\prod_{j=0}^{k} \sin(\pi(1-s+\alpha j))}.
\end{equation}
Now we need to prove that as $N \to+\infty$ the integral in the right-hand
side of (\ref{p_shift_contour2}) converges to zero for all $x \in\rr
^+$. Intuitively this is clear, since the ratio of gamma functions
in (\ref{mm_Ckl}) decreases to zero faster than any exponential
function as $\re(s) \to-\infty$, while the other factor is just
a ratio of periodic functions in $s$. The rigorous proof can be
obtained as follows:
\[
\biggl| \int_{c_0-N+\ii\rr} \mm(s) x^{-s} \,\dd s \biggr|<
C x^{N-c_0} \int_{\rr} \biggl| \frac{\Gamma(({N-c_0+\ii
t})/{\alpha} )}
{\Gamma(N-c_0+\ii t)} \biggr| g(t) \,\dd t,
\]
where we have denoted
\[
g(t)= e^{{\pi}/{\alpha} l |t|}
\prod_{j=0}^{k} \bigl| \operatorname{cosech} \bigl(\pi\bigl(t+\ii(\alpha
j-c_0)\bigr)\bigr) \bigr|.
\]
Using Stirling's approximation for the gamma function one can check
that for all $x>0$, the function
\[
x^{N-1} \frac{\Gamma(({N-c_0+\ii t})/{\alpha} )}
{\Gamma(N-c_0+\ii t)}
\]
converges to zero as $N \to+\infty$ (uniformly in $t \in\rr$); thus
the integral in the right-hand side of (\ref{p_shift_contour2})
vanishes as $N \to+\infty$, and we have a convergent series
representation (\ref{alpha_12_eq_1}). The convergence of series (\ref
{alpha_01_eq_1}) can be established in the same way, except that now we
have to shift the contour of integration in the opposite direction.
\end{pf*}

It is important to note that all the asymptotic expansions and series
representations for $p(x)$ presented in Theorems \ref{thm1} and \ref{thm2}
can be differentiated $N$ times term-by-term, where $N\ge1$ is an
arbitrary integer.
For the series representation this follows easily by the standard
argument of interchanging derivative and summation (both the series and
its derivatives converge
uniformly on compact subsets of $\rr^+$).
More work is needed to establish a similar result for the
asymptotic expansions, as it is not generally true that one can
differentiate asymptotic expansions term-by-term. A classic
counter-example is provided by the function
\[
f(x)=\frac{1}{1-x}+ e^{-x} \cos(e^{2x}) \sim\sum
_{n\ge0} x^{-n}\qquad \mbox{as } x\to+\infty,
\]
for which the asymtptotic expansion of $f'(x)$ cannot be obtained by
simply taking term-by-term derivative of the asymptotic expansion of $f(x)$.
Thus in order to establish the result on the asymptotic expansion of
$p^{(N)}(x)$ we would have to repeat the steps of the proof of Theorem
\ref{thm1}: first we take $N$th derivative of both sides of (\ref
{p_as_inv_of_m}). Then we interchange the order of integration and
differentiation in the right-hand side of (\ref{p_as_inv_of_m}); this
can be easily justified using the standard uniform-convergence
argument. Finally we shift the contour of integration to obtain
asymptotic estimates.

It is very likely that the asymptotic expansions given in Theorem \ref
{thm1} can also be interpreted as convergent series;
however, it seems to be very hard to prove this fact analytically. The
above proof of Theorem \ref{thm2} was based
on two facts: $\mm(s)$ decays faster than any exponential function
as $\re(s) \to\pm\infty$ [the sign depends on whether $\alpha\in
(0,1)$ or $\alpha\in(1,2)$], and that $\mm(s)$ is essentially a
product of
a function which decays very fast and a function which is periodic. The
first fact is still true in the general case: using (\ref{G_asympt_formula})
and (\ref{formula_Ms_general}) one can prove that
as $\re(s)\to\infty$ and $0<\epsilon<|{\arg}(s)|<\pi- \epsilon$,
we have
\begin{eqnarray*}
\log(\mm(s)) &=& (s-1) \ln(s) \biggl( 1-\frac{1}{\alpha}\biggr) +O(s), \\
\log(\mm(-s)) &=& -(s+1) \ln(s) \biggl( 1-\frac{1}{\alpha}\biggr) +O(s).
\end{eqnarray*}
The major problem now is that when $X_t \notin\cc_{k,l}$ we do not
have any periodicity. Moreover, to make matters worse,
when we move the contour of integration farther away from zero, the
poles of $\mm(s)$ become more and more dense,
and it is very hard to find an upper bound on $|\mm(s)|$ for small
values of $\im(s)$.

One can also see that the behavior of the coefficients $a_{m,n}$ and
$b_{m,n}$ is much more unpredictable compared to $c_{m,n}$, and it is
hard to say anything about the growth/decay of these coefficients as
$m$ or $n$ becomes large. Numerical results, however, indicate that the
product of ratios of $\sin(\cdot)$ functions in (\ref{def_a_mn})
remains bounded as $m$ or $n$ becomes large,
as long as $\alpha$ is not too close to a rational number. Thus it
seems reasonable to expect that the following conjecture is true:

\textit{Conjecture}: Assume that $\alpha\notin{\mathcal L} \cup\q$. If
$\alpha\in(1,2)$ $\{\alpha\in(0,1)\}$, then the
infinite series (\ref{eqn_p_0}) \{(\ref{eqn_p_infty})\} converges to
$p(x)$ for all $x>0$.

\section*{Acknowledgments}
The author would like to thank two anonymous referees
for many detailed comments and constructive suggestions.

\begin{supplement}
\stitle{Appendix A: Detailed proofs of some results related to the double gamma function}
\slink[doi]{10.1214/10-AOP577SUPP}
\sdatatype{.pdf}
\sfilename{aop577\_suppl.pdf}
\sdescription{This supplement material provides detailed computations needed to derive formulas
(4.10), (4.11), (7.1), (7.2), (7.5)
and to prove Corollary 3 and Theorem 8.}
\end{supplement}

%

%
\printaddresses


\begin{thebibliography}{26}

\bibitem{Andrews}
%
\begin{bbook}[mr]
\bauthor{\bsnm{Andrews},~\bfnm{George~E.}\binits{G.~E.}},
\bauthor{\bsnm{Askey},~\bfnm{Richard}\binits{R.}} \AND
\bauthor{\bsnm{Roy},~\bfnm{Ranjan}\binits{R.}}
(\byear{1999}).
\btitle{Special Functions}.
\bseries{Encyclopedia of Mathematics and Its Applications}
\bvolume{71}.
\bpublisher{Cambridge Univ. Press}, \baddress{Cambridge}.
\bid{mr={1688958}}
\end{bbook}
%
\endbibitem

\bibitem{Barnes1899}
%
\begin{barticle}[vtex]
\bauthor{\bsnm{Barnes},~\bfnm{E.~W.}\binits{E.~W.}}
(\byear{1899}).
\btitle{The genesis of the double gamma function}.
\bjournal{Proc. London Math. Soc.}
\bvolume{31}
\bpages{358--381}.
\end{barticle}
%
\endbibitem

\bibitem{Barnes1901}
%
\begin{barticle}[vtex]
\bauthor{\bsnm{Barnes},~\bfnm{E.~W.}\binits{E.~W.}}
(\byear{1901}).
\btitle{The theory of the double gamma function}.
\bjournal{Phil. Trans. Royal Soc. London (A)}
\bvolume{196}
\bpages{265--387}.
\end{barticle}
%
\endbibitem

\bibitem{Bernyk2008}
%
\begin{barticle}[mr]
\bauthor{\bsnm{Bernyk},~\bfnm{Violetta}\binits{V.}},
\bauthor{\bsnm{Dalang},~\bfnm{Robert~C.}\binits{R.~C.}} \AND
\bauthor{\bsnm{Peskir},~\bfnm{Goran}\binits{G.}}
(\byear{2008}).
\btitle{The law of the supremum of a stable {L}\'evy process with no negative
jumps}.
\bjournal{Ann. Probab.}
\bvolume{36}
\bpages{1777--1789}.
\bid{doi={10.1214/07-AOP376}, mr={2440923}}
\end{barticle}
%
\endbibitem

\bibitem{Bertoin}
%
\begin{bbook}[mr]
\bauthor{\bsnm{Bertoin},~\bfnm{Jean}\binits{J.}}
(\byear{1996}).
\btitle{L\'evy Processes}.
\bseries{Cambridge Tracts in Mathematics}
\bvolume{121}.
\bpublisher{Cambridge Univ. Press}, \baddress{Cambridge}.
\bid{mr={1406564}}
\end{bbook}
%
\endbibitem

\bibitem{Bill1997}
%
\begin{barticle}[mr]
\bauthor{\bsnm{Billingham},~\bfnm{J.}\binits{J.}} \AND
\bauthor{\bsnm{King},~\bfnm{A.~C.}\binits{A.~C.}}
(\byear{1997}).
\btitle{Uniform asymptotic expansions for the {B}arnes double gamma function}.
\bjournal{Proc. Roy. Soc. London Ser. A}
\bvolume{453}
\bpages{1817--1829}.
\bid{doi={10.1098/rspa.1997.0098}, mr={1478136}}
\end{barticle}
%
\endbibitem

\bibitem{Bingham1975}
%
\begin{barticle}[vtex]
\bauthor{\bsnm{Bingham},~\bfnm{N.~H.}\binits{N.~H.}}
(\byear{1975}).
\btitle{Fluctuation theory in continuous time}.
\bjournal{Adv. in Appl. Probab.}
\bvolume{7}
\bpages{705--766}.
\bid{mr={0386027}}
\end{barticle}
%
\endbibitem

\bibitem{Borodin}
%
\begin{bbook}[vtex]
\bauthor{\bsnm{Borodin},~\bfnm{Andrei~N.}\binits{A.~N.}} \AND
\bauthor{\bsnm{Salminen},~\bfnm{Paavo}\binits{P.}}
(\byear{1996}).
\btitle{Handbook of {B}rownian Motion---Facts and Formulae}.
\bpublisher{Birkh\"auser}, \baddress{Basel}.
\bid{mr={1477407}}
\end{bbook}
%
\endbibitem

\bibitem{Borwein2000}
%
\begin{barticle}[vtex]
\bauthor{\bsnm{Borwein},~\bfnm{Jonathan~M.}\binits{J.~M.}},
\bauthor{\bsnm{Bradley},~\bfnm{David~M.}\binits{D.~M.}} \AND
\bauthor{\bsnm{Crandall},~\bfnm{Richard~E.}\binits{R.~E.}}
(\byear{2000}).
\btitle{Computational strategies for the {R}iemann zeta function}.
\bjournal{J. Comput. Appl. Math.}
\bvolume{121}
\bpages{247--296}.
\bid{doi={10.1016/S0377-0427(00)00336-8}, mr={1780051}}
\end{barticle}
%
\endbibitem

\bibitem{Darling1956}
%
\begin{barticle}[mr]
\bauthor{\bsnm{Darling},~\bfnm{D.~A.}\binits{D.~A.}}
(\byear{1956}).
\btitle{The maximum of sums of stable random variables}.
\bjournal{Trans. Amer. Math. Soc.}
\bvolume{83}
\bpages{164--169}.
\bid{mr={0080393}}
\end{barticle}
%
\endbibitem

\bibitem{Doney1987}
%
\begin{barticle}[vtex]
\bauthor{\bsnm{Doney},~\bfnm{R.~A.}\binits{R.~A.}}
(\byear{1987}).
\btitle{On {W}iener--{H}opf factorisation and the distribution of
extrema for
certain stable processes}.
\bjournal{Ann. Probab.}
\bvolume{15}
\bpages{1352--1362}.
\bid{mr={905336}}
\end{barticle}
%
\endbibitem

\bibitem{Doney2008}
%
\begin{barticle}[mr]
\bauthor{\bsnm{Doney},~\bfnm{R.~A.}\binits{R.~A.}}
(\byear{2008}).
\btitle{A note on the supremum of a stable process}.
\bjournal{Stochastics}
\bvolume{80}
\bpages{151--155}.
\bid{doi={10.1080/17442500701830399}, mr={2402160}}
\end{barticle}
%
\endbibitem

\bibitem{Doney2010}
%
\begin{barticle}[mr]
\bauthor{\bsnm{Doney},~\bfnm{R.~A.}\binits{R.~A.}} \AND
\bauthor{\bsnm{Savov},~\bfnm{M.~S.}\binits{M.~S.}}
(\byear{2010}).
\btitle{The asymptotic behavior of densities related to the supremum
of a
stable process}.
\bjournal{Ann. Probab.}
\bvolume{38}
\bpages{316--326}.
\bid{doi={10.1214/09-AOP479}, mr={2599201}}
\end{barticle}
%
\endbibitem

\bibitem{Fourati}
%
\begin{barticle}[mr]
\bauthor{\bsnm{Fourati},~\bfnm{S.}\binits{S.}}
(\byear{2006}).
\btitle{Inversion de l'espace et du temps des processus de {L}\'evy stables}.
\bjournal{Probab. Theory Related Fields}
\bvolume{135}
\bpages{201--215}.
\bid{doi={10.1007/s00440-005-0455-2}, mr={2218871}}
\end{barticle}
%
\endbibitem

\bibitem{Graczyk2009}
%
\begin{barticle}[vtex]
\bauthor{\bsnm{Graczyk},~\bfnm{P.}\binits{P.}} \AND
\bauthor{\bsnm{Jakubowski},~\bfnm{T.}\binits{T.}}
(\byear{2009}).
\btitle{Wiener--Hopf factors for stable processes}.
\bjournal{Ann. Inst. H. Poincar\'e Probab. Statist.}
\bnote{To appear}.
\end{barticle}
%
\endbibitem

\bibitem{Heyde1969}
%
\begin{barticle}[vtex]
\bauthor{\bsnm{Heyde},~\bfnm{C.~C.}\binits{C.~C.}}
(\byear{1969}).
\btitle{On the maximum of sums of random variables and the supremum functional
for stable processes}.
\bjournal{J. Appl. Probab.}
\bvolume{6}
\bpages{419--429}.
\bid{mr={0251766}}
\end{barticle}
%
\endbibitem

\bibitem{Jeffrey2007}
%
\begin{bbook}[vtex]
\bauthor{\bsnm{Jeffrey},~\bfnm{A.}\binits{A.}~ed.}
(\byear{2007}).
\btitle{Table of Integrals, Series and Products}, \bedition{7th} ed.
\bpublisher{Academic Press}, \baddress{Amsterdam}.
\end{bbook}
%
\endbibitem

\bibitem{Khinchin}
%
\begin{bbook}[vtex]
\bauthor{\bsnm{Khinchin},~\bfnm{A.~Ya.}\binits{A.~Y.}}
(\byear{1997}).
\btitle{Continued Fractions}, \bedition{Russian} ed.
\bpublisher{Dover Publications Inc.}, \baddress{Mineola, NY}.
\bid{mr={1451873}}
\end{bbook}
%
\endbibitem

\bibitem{Kuznetsov}
%
\begin{barticle}[vtex]
\bauthor{\bsnm{Kuznetsov},~\bfnm{A.}\binits{A.}}
(\byear{2009}).
\btitle{Analytical proof of
Pecherskii--Rogozin identity and Wiener--Hopf factorization}.
\bjournal{Theory Probab. Appl.}
\bnote{To appear}.
\end{barticle}
%
\endbibitem

\bibitem{Lawrie1994}
%
\begin{barticle}[mr]
\bauthor{\bsnm{Lawrie},~\bfnm{J.~B.}\binits{J.~B.}} \AND
\bauthor{\bsnm{King},~\bfnm{A.~C.}\binits{A.~C.}}
(\byear{1994}).
\btitle{Exact solution to a class of functional difference equations with
application to a moving contact line flow}.
\bjournal{European J. Appl. Math.}
\bvolume{5}
\bpages{141--157}.
\bid{doi={10.1017/S0956792500001364}, mr={1285035}}
\end{barticle}
%
\endbibitem

\bibitem{Lewin1981}
%
\begin{bbook}[vtex]
\bauthor{\bsnm{Lewin},~\bfnm{Leonard}\binits{L.}}
(\byear{1981}).
\btitle{Polylogarithms and Associated Functions}.
\bpublisher{North-Holland}, \baddress{New York}.
\bid{mr={618278}}
\end{bbook}
%
\endbibitem

\bibitem{Mordecki}
%
\begin{barticle}[vtex]
\bauthor{\bsnm{Lewis},~\bfnm{Alan~L.}\binits{A.~L.}} \AND
\bauthor{\bsnm{Mordecki},~\bfnm{Ernesto}\binits{E.}}
(\byear{2008}).
\btitle{Wiener--{H}opf factorization for {L}\'evy processes having positive
jumps with rational transforms}.
\bjournal{J. Appl. Probab.}
\bvolume{45}
\bpages{118--134}.
\bid{doi={10.1239/jap/1208358956}, mr={2409315}}%
\end{barticle}%
%
\endbibitem%

\bibitem{Patie2009}
%
\begin{barticle}[mr]
\bauthor{\bsnm{Patie},~\bfnm{P.}\binits{P.}}
(\byear{2009}).
\btitle{A few remarks on the supremum of stable processes}.
\bjournal{Statist. Probab. Lett.}
\bvolume{79}
\bpages{1125--1128}.
\bid{doi={10.1016/j.spl.2009.01.001}, mr={2510779}}
\end{barticle}
%
\endbibitem

\bibitem{Peskir2008}
%
\begin{barticle}[vtex]
\bauthor{\bsnm{Peskir},~\bfnm{Goran}\binits{G.}}
(\byear{2008}).
\btitle{The law of the hitting times to points by a stable {L}\'evy process
with no negative jumps}.
\bjournal{Electron. Comm. Probab.}
\bvolume{13}
\bpages{653--659}.
\bid{mr={2466193}}
\end{barticle}
%
\endbibitem

\bibitem{Vatutin2009}
%
\begin{barticle}[mr]
\bauthor{\bsnm{Vatutin},~\bfnm{Vladimir~A.}\binits{V.~A.}} \AND
\bauthor{\bsnm{Wachtel},~\bfnm{Vitali}\binits{V.}}
(\byear{2009}).
\btitle{Local probabilities for random walks conditioned to stay positive}.
\bjournal{Probab. Theory Related Fields}
\bvolume{143}
\bpages{177--217}.
\bid{doi={10.1007/s00440-007-0124-8}, mr={2449127}}
\end{barticle}
%
\endbibitem

\bibitem{Zolotarev1957}
%
\begin{barticle}[vtex]
\bauthor{\bsnm{Zolotarev},~\bfnm{V.~M.}\binits{V.~M.}}
(\byear{1957}).
\btitle{Mellin--{S}tieltjes transformations in probability theory}.
\bjournal{Teor. Veroyatnost. i Primenen.}
\bvolume{2}
\bpages{444--469}.
\bid{mr={0108843}}
\end{barticle}
%
\endbibitem

\end{thebibliography}
\end{document}